\documentclass[12pt,reqno]{article}
\usepackage[pdftex]{hyperref}
\usepackage{amsmath, amsthm, mathrsfs, graphicx,amsfonts, amssymb,color}
\newtheorem{thm}{Theorem}[section]
\newtheorem{cor}[thm]{Corollary}
\newtheorem{lem}[thm]{Lemma}
\newtheorem{prop}[thm]{Proposition}
\theoremstyle{Definition}
\newtheorem{defn}[thm]{Definition}

\theoremstyle{example}
\newtheorem{exm}[thm]{Example}
\numberwithin{equation}{section}

\def\dsum{\displaystyle\sum}
\def\dsup{\displaystyle\sup}
\def\dlim{\displaystyle\lim}
\def\dlimsup{\displaystyle\limsup}
\def\dmax{\displaystyle\max}
\def\dmin{\displaystyle\min}
\def\dinf{\displaystyle\inf}

\newcommand{\scr}[1]{\mathscr #1}
\newcommand{\norm}[2]{\left\|{#1}\right\|_{#2}}
\newcommand{\abs}[1]{\left\vert#1\right\vert}
\newcommand{\set}[1]{\left\{#1\right\}}
\newcommand{\R}{\mathbb R}
\newcommand{\eps}{\varepsilon}
\newcommand{\A}{\mathcal{A}}
\newcommand{\E}{\mathbb{E}}
\newcommand{\D}{\scr{D}}
\renewcommand{\P}{\mathbb P}
\newcommand{\nnb}{\nonumber}

\def\R{\mathbb R}
\def\P{\mathbb P}
\def\Z{\mathbb Z}
\def\N{\mathbb N}
\def\F{\scr F}
\def\K{\scr K}
\def\bg{\begin}
\def\be{\bg{equation}}
\def\de{\end{equation}}
\def\bgar{\bg{eqnarray}}
\def\edar{\end{eqnarray}}
\def\beqnn{\begin{eqnarray*}}
\def\eeqnn{\end{eqnarray*}}
\def\lb{\label}
\def\ct{\cite}
\def\l{\left}
\def\r{\right}
\def\fr{\frac}
\def\alp{\alpha}
\def\bt{\beta}
\def\gm{\gamma}
\def\Gm{\Gamma}
\def\dlt{\delta}
\def\Dlt{\Delta}
\def\eps{\epsilon}
\def\veps{\varepsilon}
\def\tht{\theta}
\def\Tht{\Theta}
\def\kp{\kappa}
\def\lmd{\lambda}
\def\Lmd{\Lambda}
\def\vro{\varrho}
\def\sgm{\sigma}
\def\Sgm{\Sigma}
\def\vph{\varphi}
\def\omg{\omega}
\def\Omg{\Omega}
\def\fa{\forall}
\def\emp{\emptyset}
\def\ex{\exists}
\def\nbl{\nabla}
\def\pat{\partial}
\def\ift{\infty}
\def\bca{\bigcap}
\def\bcu{\bigcup}
\def\lar{\leftarrow}
\def\Lar{\Leftarrow}
\def\rar{\rightarrow}
\def\Rar{\Rightarrow}
\def\lla{\longleftarrow}
\def\Lla{\Longleftarrow}
\def\to{\longrightarrow}
\def\To{\Longrightarrow}
\def\lra{\leftrightarrow}
\def\Lra{\Leftrightarrow}
\def\llra{\longleftrightarrow}
\def\Llra{\Longleftrightarrow}
\def\q{\quad}
\def\gap{\text {\rm gap}}
\def\var{\text {\rm Var}}
\def\V{\text {\rm V}}
\def\TV{\text {\rm TV}}
\def\ess{{\rm ess}}
\def\hess{{\rm Hess}}
\def\ric{{\rm Ric}}
\def\tr{{\rm tr}}
\def\d{{\mbox{\rm d}}}\def\e{{\mbox{\rm e}}}
\def\supp{{\mbox{\rm supp}}}
\def\lan{\langle}
\def\ran{\rangle}
\def\[{\l[} \def\]{\r]}
\def\({\l(} \def\){\r)}
\def\|{\bigg|}
\def\hat{\widehat}
\def\bar{\overline}
\def\tld{\widetilde}
\def\mpb{\vskip6pt}

\renewcommand{\aa}[3]{{#1}_{#2}({#3})}
\newcommand{\p}[2]{p_{#1}({#2})}
\newcommand{\h}[2]{h_{#1}({#2})}
\newcommand{\m}[2]{m_{#1}^{(#2)}}
\newcommand{\mb}[1]{\fr{1}{\mu_{#1}b_{#1}}}
\newcommand{\qq}[1]{q_{#1}}
\renewcommand{\d}[2]{d_{#1}^{(#2)}}
\newcommand{\he}[2]{\sum_{#1}^{#2}}
\newcommand{\x}[2]{x_{#1}^{(#2)}}
\newcommand{\rf}[1]{(\ref{#1})}
\newcommand{\pfthm}[1]{\vskip.5cm \noindent\emph{Proof of Theorem \ref{#1}}}
\newcommand{\pfcor}[1]{\vskip.5cm \noindent\emph{Proof of Corollary \ref{#1}}}
\newcommand{\pfprop}[1]{\vskip.5cm \noindent\emph{Proof of Proposition \ref{#1}}}

\newcommand{\red}[1]{{\color{red} #1}}
\newcommand{\yel}[1]{{\color{yellow} #1}}
\newcommand{\blu}[1]{{\color{blue} #1}}
\newcommand{\grn}[1]{{\color{green} #1}}

\title{{\bf  On geometric and algebraic transience for discrete-time Markov chains}}

\author{
{\bf Yong-Hua Mao and Yan-Hong Song\footnote{Correspondence should be addressed to Yan-Hong Song
(email: songyh@mail.bnu.edu.cn)}}\\
\footnotesize{School of Mathematical Sciences, Beijing Normal University, }\\
\footnotesize{Laboratory of Mathematics and Complex Systems, Ministry of Education}\\
 \footnotesize{Beijing 100875, China}\\
\footnotesize{Email:
maoyh@bnu.edu.cn, songyh@mail.bnu.edu.cn}
}
\date{ }

\begin{document}

\maketitle

\begin{abstract}
General characterizations of ergodic Markov chains
have been developed in considerable detail. In this paper, we study the transience
for discrete-time Markov chains on general state spaces, including the geometric transience and algebraic transience.
Criteria are presented through establishing the drift condition and considering
the first return time. As an application, we give explicit criteria for the random
walk on the half line and the skip-free chain on nonnegative integers.
\end{abstract}

\noindent


\section{Introduction}
In the past decades, great efforts have been made to study the ergodic
theory for Markov chains.
The drift condition (Foster-Lyapunov condition) is an important method,
which has been used extensively. For example, Meyn and Tweedie \ct{m} gave drift conditions for geometric and uniform ergodicity.
Tuominen and Tweedie \ct{tt} studied subgeometric ergodicity by using a sequence of drift conditions, which is a foundational work. Building on it, Jarner and Roberts \ct{jr} investigated polynomial ergodicity by establishing a single drift condition, and Mao \ct{m1, m3} used one drift condition to study the algebraic convergence and
the ergodic degree.
Then Douc, Fort, Moulines and Soulier \ct{dfm} presented a new practical drift condition to prove subgeometric ergodicity. This condition, extending the condition introduced by Jarner and Roberts,
turned out to be more convenient than that in Tuominen and Tweedie \ct{tt}.

In this paper, we aim to investigate the transient theory for discrete-time
Markov chains, which is also an interesting and challenging problem.
The study of transient theory may be dated back to Harris \ct{har} in the 1950s,
who obtained a necessary condition and a sufficient condition for the existence
of stationary measures for transient Markov chains.
The problem was further discussed by Vere-Jones \ct{vj, vj1},
who defined the geometric transience on the countable state space, and
studied the $\lmd$-subinvariant measure of geometrically transient Markov chains.
For more details, one can refer to \ct{and}.
Then Vere-Jones's results were extended to chains with fixed absorbing points, see \ct{fla} and references within (if the absorption is reducible or not certain, see e.g. \ct{vdo, ppk}). In \ct{tw, twee}, Tweedie extended the results of Harris and Vere-Jones to the general state space.
Based on these works, Meyn and Tweedie \ct{meyn} systematically studied
the stochastic stability of discrete-time Markov chains. In their book,
they used the drift condition to study the criteria of transience,
see \ct[Theorem 8.0.2]{meyn}.
Besides, the transient theory has a wide range of applications,
see e.g. \ct{vdf, vds, cay}.

However, in spite of these developments in both the drift condition of ergodicity
and the transient theory,
it seems that using drift conditions to study further transience of discrete-time Markov chains
has not been fully revealed. The goal of this paper is therefore to study
the geometric transience and algebraic transience (see Definitions
\ref{fgh} and \ref{fhg} below) of general discrete-time
Markov chains, through establishing appropriate drift conditions.

Let us introduce the basic setup of the paper.
Let $\Phi=\{\Phi_n: n\in\Z_+\}$ be a discrete-time
homogeneous Markov chain on a general state space $X$,
endowed with a countably generated $\sigma$-field $\mathcal{B}(X)$.
Denote by $P^{n}(x, A)$
the $n$-step transition kernel of the chain:
$$P^{n}(x, A)=\P_{x}\{\Phi_{n}\in A\},\q n\in \Z_+,\; x\in X,\; A\in\mathcal{B}(X),$$
where $\P_x$ is the conditional distribution of the chain given $\Phi_0=x$.
The corresponding expectation operator will be denoted $\E_x$.
Here, $P$ may be stochastic or sub-stochastic, and
for all nonnegative measurable function $f$,
$$P^{n}f(x)=\int_X f(y)P^{n}(x, d y),\q n\in\Z_+, \; x\in X.$$
Assume throughout the paper
that the chain $\Phi$ is $\psi$-irreducible, where $\psi$ is a maximal
irreducibility measure.
Write $\mathcal{B}^{+}(X)=\l\{A\in\mathcal{B}(X): \psi(A)>0\r\}$
for the sets of positive $\psi$-measure.

For a probability distribution $a=(a_n)_{n\in\N}$,
let $K_a$ be the transition kernel given by
$$K_a(x, A)=\sum_{n=1}^{\ift}a_n P^{n}(x, A),\q x\in X,\; A\in\mathcal{B}(X).$$
A set $A\in\mathcal{B}(X)$ is called petite if
there exists a probability distribution $a$ and a nontrivial measure $\nu_a$ such that
$$K_a(x, \cdot)\geq\nu_a(\cdot),\q x\in A.$$
Petite sets are not rare:
if $\Phi$ is $\psi$-irreducible, then for every $B\in\mathcal{B}^{+}(X)$,
there exists a petite set $A\subset B$ such that $A\in\mathcal{B}^{+}(X)$, see
\ct[Theorem 5.2.2]{meyn} for reference.

The first return time of a set $A\in \mathcal{B}(X)$ is denoted by
$\tau_{A}=\inf\{n\geq1: \Phi_{n}\in A\}$,
and the first hitting time is defined by $\sigma_{A}=\tau_A 1_{\l\{\Phi_0\notin A\r\}}=\inf\{n\geq0: \Phi_{n}\in A\}$.
They are two stopping times with respect to the filtration $(\mathcal{F}_n)$,
where $\mathcal{F}_n=\sgm\l\{\Phi_0, \cdots, \Phi_n\r\}$. Let
\be\lb{j8}
F^{n}(x, A)=\P_x\{\tau_{A}=n\},\q n=1, 2, \cdots, \ift\nnb
\de
be the distribution of $\tau_A$, and
$$L(x, A)=\sum_{n=1}^{\ift}F^{n}(x, A)=\P_x\{\tau_A<\ift\}$$
the probability of $\Phi$ ever returning to $A$.

Recall that the chain $\Phi$ is transient if it is $\psi$-irreducible and there exist sets $A_i\in\mathcal{B}^{+}(X)$, $i=1$, $2$, $\cdots$ such that
\bg{equation}\lb{cr68}
X=\bigcup_{i=1}^{\ift}A_i\q\mbox{and}\q\sup_{x\in A_i}\sum_{n=1}^{\ift}P^{n}(x, A_i)<\ift,\q i\geq 1.
\end{equation}
Moreover, according to the proof of \ct[Theorem 8.3.6]{meyn}, we can have
\bg{prop}\lb{ti8}
The chain $\Phi$ is transient if and only if for every
petite set $B\in\mathcal{B}^{+}(X)$, there exists a set $A\subset B$ with
$\psi(A)>0$ such that
\bg{equation}\lb{cris0}
\sup_{x\in A}L(x, A)<1.\nnb
\end{equation}
\end{prop}
For the transient chain, by \rf{cr68}, we have
$\lim_{n\rightarrow\ift}P^{n}(x, A_i)=0$ for all $x\in A_i$.
Thus, it is natural to ask how fast $P^{n}(x, A_i)$ goes to zero.
This is the main motivation for us to study further transience, which
we specify to be geometric transience and algebraic transience.
In the paper,
we will give practical drift conditions for these transience,
as have been done in the ergodic case. The basic idea is still to consider
the first return time. Let us take as an example the comparison
of geometric ergodicity and geometric transience.

The chain $\Phi$ is called geometrically ergodic
if there exists a stationary distribution $\pi$ satisfying
$$||P^{n}(x, \cdot)-\pi||\leq M(x)\rho^{n},\q n\in\Z_+,\; x\in X,$$
for some $M(x)<\ift$ and $\rho<1$, where $||\cdot||$ is
the total variation norm. For the ergodicity, Meyn and Tweedie \ct[Chapter 15]{meyn}
have the following main results.
\bg{thm}\lb{bta}
Suppose that the chain $\Phi$ is $\psi$-irreducible and aperiodic.
Then the following statements are equivalent.

$(1)$ There exist some petite set $A\in\mathcal{B}^{+}(X)$ and $\kp>1$ such that
\be\lb{xi0}
\sup_{x\in A}\E_x\l[\kp^{\tau_A}\r]<\ift.
\end{equation}

$(2)$ There exist some petite set $A\in\mathcal{B}^{+}(X)$, constants $b<\ift$, $\lmd<1$
and a function $W\geq 1$, with $W(x_0)<\ift$ for some $x_0\in X$, satisfying
the drift condition
\be\lb{si81}
P W(x)\leq\lmd W(x)+b 1_A(x),\q x\in X\nnb.
\end{equation}

$(3)$ The chain $\Phi$ is geometrically ergodic.
\end{thm}
Note that
$L(x, A)=1$ for the ergodic Markov chain, we can rewrite \rf{xi0} as
\be\lb{l9a}
\sup_{x\in A}L(x, A)=1 \q\mbox{and}\q \sup_{x\in A}\E_x\l[\kp^{\tau_A}1_{\{\tau_A<\ift\}}\r]<\ift.
\end{equation}
As for the geometric transience,
it shows in the following Theorem \ref{zzz} that if (and only if)
\be\lb{k9a}
\sup_{x\in A}L(x, A)<1 \q\mbox{and}\q \sup_{x\in A}\E_x\l[\kp^{\tau_A}1_{\{\tau_A<\ift\}}\r]<\ift,
\end{equation}
then the chain $\Phi$ is geometrically transient.
Thus, from \rf{l9a} and \rf{k9a}, it is natural to study $F^{n}(x, A)$ more carefully for the geometric transience.

The remainder of the paper is organized as follows.
The geometric transience, including strongly geometric transience and
uniformly geometric transience are investigated in Section \ref{geo}.
Section \ref{als} is devoted to researching the algebraic transience.
In Section \ref{exan}, we apply our results to the random walk on $\R_+$
and the skip-free chain on $\Z_+$.


\section{Geometric transience}\label{geo}

In this section, we will study three kinds of geometric transience.

\subsection{Geometric transience}\label{pt}

We begin with the definition of geometric transience.

\bg{defn}\lb{fgh}
A set $A\in \mathcal{B}^{+}(X)$ is called uniformly geometric transient
if there exists a constant $\kappa>1$ such that
\be\lb{vfr}
\sup_{x\in A}\sum_{n=1}^{\infty}\kappa^{n}P^{n}(x, A)<\ift\nnb.
\de
The chain $\Phi$ is called geometrically transient
if it is $\psi$-irreducible and $X$ can be covered $\psi$-a.e.
by a countable number of uniformly geometric transient sets.
That is, there exist sets $D$ and $A_i$, $i=1, 2, \cdots$ such that
$X=D\cup\l(\bigcup_{i=1}^{\ift}A_i\r)$, where $\psi(D)=0$ and each $A_i$ is
uniformly geometric transient.
\end{defn}

For the geometric transience, we have the following main result linking
the ``local" geometric transience, the first return time,
the drift condition and the geometric transience.
\bg{thm}\lb{zzz}
Suppose that the chain $\Phi$ is $\psi$-irreducible.
Then the following statements are equivalent.

$(1)$ There exist some set $A\in\mathcal{B}^{+}(X)$ and $\kp>1$ such that
$$\sup_{x\in A}\sum_{n=1}^{\ift}\kp^{n}P^{n}(x, A)<\ift.$$

$(2)$ There exist some set $A\in\mathcal{B}^{+}(X)$ and $\kp>1$ such that
\be\lb{xi01}
\sup_{x\in A}\E_x\l[\kp^{\tau_A}1_{\{\tau_A<\ift\}}\r]<1.
\end{equation}

$(3)$ There exist some set $A\in\mathcal{B}^{+}(X)$ and $\kp>1$ such that
$$\sup_{x\in A}L(x, A)<1,\q \sup_{x\in A}\E_x\l[\kp^{\tau_A}1_{\{\tau_A<\ift\}}\r]<\ift.$$

$(4)$ There exist some set $A\in\mathcal{B}^{+}(X)$, constants $b$, $\lmd\in(0, 1)$,
and a function $W\geq 1_A$,
with $W(x_0)<\ift$ for some $x_0\in X$, satisfying the drift condition
\be\lb{34e}
P W(x)\leq\lmd W(x) 1_{A^{c}}(x)+b 1_A(x),\q x\in X.
\end{equation}

$(5)$ The chain $\Phi$ is geometrically transient.

\end{thm}

\bg{rem}\lb{s683}
$(1)$ According to the proof of $(5)\Rar(3)$, the set $A\in\mathcal{B}^{+}(X)$ is a petite set.

$(2)$ Since $P W(x)\leq b$ holds for $x\in A$ with $b\in(0, 1)$, the set $\{x\in A^{c}: W(x)<1\}\not=\emptyset$ when $P$ is stochastic.
\end{rem}

In order to prove the theorem, we need three lemmas.
Let $\Lambda$ be the family of increasing functions
$r$: $\mathbb{Z}_{+}\rightarrow [1, \infty)$ satisfying
\be\lb{lk}
r(0)=1
\q\mbox{and}\q r(m+n)\leq r(m)r(n), \q m,\, n\in \mathbb{Z}_{+}\nnb.
\de
The next lemma is a straightforward generalization of \ct[Proposition 2.1]{tw}.

\begin{lem}\lb{nc}
Let $r\in\Lmd$.

$(1)$
Assume that there exists a set
$A\in \mathcal{B}^{+}(X)$ such that
$$\sum_{n=1}^{\ift}r(n)P^{n}(x, A)<\ift,\q x\in A.$$
Then there exist sets $D$ and $A_i$, $i=1, 2, \cdots$
such that $X=D\cup\l(\bigcup_{i=1}^{\ift}A_i\r)$, $\psi(D)=0$, and
$$\sup_{x\in A_i}\sum_{n=1}^{\ift}r(n)P^{n}(x, A_i)<\ift,\q i\geq1.$$

$(2)$
Assume that there exists a set
$A\in \mathcal{B}^{+}(X)$ such that
$$\sup_{x\in X}\sum_{n=1}^{\ift}r(n)P^{n}(x, A)<\ift.$$
Then there exist sets $A_i$, $i=1, 2, \cdots$
such that $X=\bigcup_{i=1}^{\ift}A_i$, and
$$\sup_{x\in X}\sum_{n=1}^{\ift}r(n)P^{n}(x, A_i)<\ift,\q i\geq1.$$
\end{lem}

\bg{proof}
We only prove the first assertion, since the proof of the second one is similar.

(a) Set
$D=\l\{x\in X: \sum_{n=1}^{\ift}r(n)P^{n}(x, A)=\ift\r\}$.
Since $r$ is increasing,
$$r(m+n)P^{m+n}(x, A)\geq\int_{D}P^{m}(x, d y)r(n)P^{n}(y, A),\q m, n\in\N.$$
Summing over $n$ gives
$$\ift>\sum_{n=1}^{\ift}r(n)P^{n}(x, A)\geq\int_{D}P^{m}(x, d y)
\sum_{n=1}^{\ift}r(n)P^{n}(y, A),\q x\in A,$$
which means $P^{m}(x, D)=0$ for $m\in\N$. Then $\psi(D)=0$
by the $\psi$-irreducibility.

(b) For $n$, $j\in\N$, set
$$H(n, j)=\l\{x\in D^{c}: P^{n}(x, A)\in \l((j+1)^{-1}, j^{-1}\r],
P^{k}(x, A)=0, k=1, 2, \cdots, n-1\r\}.$$
Then
$D^{c}=\bigcup_{n, j=1}^{\ift}H(n, j)$ by the $\psi$-irreducibility.
Using $r(m+n)\geq r(n)$ again, we have
$$\aligned
r(m+n)P^{m+n}(x, A)&\geq \int_{H(n, j)}r(m)P^{n}(y, A)P^{m}(x, d y)\geq (j+1)^{-1}r(m)P^{m}(x, H(n, j)).\\
\endaligned
$$
Summing over $m$ gives $$\sum_{m=1}^{\ift}r(m)P^{m}(x, A)
\geq (j+1)^{-1}\sum_{m=1}^{\ift}r(m)P^{m}(x, H(n, j)).$$
Hence $\sum_{m=1}^{\ift}r(m)P^{m}(x, H(n, j))<\ift$ for $x\in D^{c}$.

(c) For $k\in\N$, let
$$B(n, j, k)=\l\{x\in H(n, j):
\sum_{m=1}^{\ift}r(m)P^{m}(x, H(n, j))\leq k\r\}.$$
Then it is obvious that $H(n, j)=\bigcup_{k=1}^{\ift}B(n, j, k)$.
Combining this with (a) and (b), we have
$X=D\cup\l(\bigcup_{n, j, k=1}^{\ift}B(n, j, k)\r)$, $\psi(D)=0$ and
$$\sup_{x\in B(n, j, k)}\sum_{m=1}^{\ift}r(m)P^{m}\l(x, B(n, j, k)\r)<\ift,\q n,\, j,\, k\in\N,$$
which yields the desired conclusion.
\end{proof}

\bg{cor}\lb{y61}
If $\Phi$ is geometrically transient, then it is transient.
\end{cor}
\bg{proof}
Suppose that $\Phi$ is geometrically transient. Then by Definition \ref{fgh},
there exist $A\in\mathcal{B}^{+}{(X)}$ and $\kp>1$ such that
$$\sup_{x\in A}\sum_{n=1}^{\infty}\kappa^{n}P^{n}(x, A)<\ift,$$
which implies $\sup_{x\in A}\sum_{n=1}^{\infty}P^{n}(x, A)<\ift$.
Thus, according to the first entrance decomposition, for all $x\in X$,
\be\lb{qq40}\aligned
\sum_{n=1}^{\ift}P^{n}(x, A)
&=\sum_{n=1}^{\ift}F^{n}(x, A)+\sum_{n=1}^{\ift}\sum_{m=1}^{n-1}\int_{A}P^{n-m}(y, A)F^{m}(x, d y)\\&=
L(x, A)+\int_{A}\sum_{n=1}^{\ift}P^{n}(y, A)L(x, d y)\\
&\leq 1+\sup_{y\in A}\sum_{n=1}^{\ift}P^{n}(y, A).\nnb\\
\endaligned
\de
That is, $\sup_{x\in X}\sum_{n=1}^{\infty}P^{n}(x, A)<\ift$.
Hence the chain is transient from Lemma \ref{nc}(2) by letting $r(n)=1$.
\end{proof}

We next give the condition on the first return time which
ensures that a set is uniformly geometric transient.

\begin{lem}\lb{n}
Let $A\in \mathcal{B}^{+}(X)$ and $\kp\geq1$.
Suppose that there exists a constant
$\varepsilon\in(0, 1)$ such that
$$\sum_{n=1}^{\ift}\kappa^{n}F^{n}(x, A)\leq \varepsilon,\q x\in A.$$
Then we have $$\sum_{n=1}^{\ift}\kappa^{n}P^{n}(x, A)\leq\frac{\varepsilon}
{1-\varepsilon},\q x\in A.$$
\end{lem}

\bg{proof}
For $A \in \mathcal{B}^{+}(X)$, the last exit decomposition can be written as
\be\lb{le}
P^{n}(x, A)=F^{n}(x, A)+
\sum_{m=1}^{n-1}\int_{A} P^{m}(x, d y)F^{n-m}(y, A),\q n\in\N.
\de
For fixed $N\in\mathbb{N}$, multiplying by $\kappa^{n}$ in (\ref{le}) and
summing $n$ from 1 to $N$, we obtain
\be\lb{qqm9}\aligned
\sum_{n=1}^{N}\kappa^{n}P^{n}(x, A)
&=\sum_{n=1}^{N}\kappa^{n}F^{n}(x, A)+
\sum_{n=1}^{N}\sum_{m=1}^{n-1}\int_{A}\kappa^{m} P^{m}(x, d y)
\kappa^{n-m}F^{n-m}(y, A)\\
&=\sum_{n=1}^{N}\kappa^{n}F^{n}(x, A)+
\int_{A}\sum_{m=1}^{N-1}\kappa^{m}P^{m}(x, d y)\sum_{n=1}^{N-m}
\kappa^{n}F^{n}(y, A)\\
&\leq\veps+
\veps\sum_{n=1}^{N}\kappa^{n}P^{n}(x, A).\nnb\\
\endaligned
\de
That is,
$\sum_{n=1}^{N}\kappa^{n}P^{n}(x, A)\leq\fr{\veps}{1-\veps}$,
which yields the assertion by letting $N\rightarrow\ift$.
\end{proof}

To investigate the drift condition for the geometric transience, we will use
the well-known minimal nonnegative solution theory,
which is an important tool to study the recurrence and transience. For more details,
one can refer to \ct{c1, hzt}.

\begin{lem}\lb{r}
For $r\in \Lambda$, set $\hat{r}(n)=\sum_{k=0}^{n}r(k)$. Let $A\in \mathcal{B}^{+}(X)$. Then $g^{*}(x):=\E_x[\hat{r}(\tau_A)1_{\{\tau_A<\ift\}}]$
is the minimal nonnegative solution of the equation
\be\lb{lsf}
g(x)= \int_{A^{c}} g(y)P(x, d y)+P(x, A)+\E_x[r(\tau_A)1_{\{\tau_A<\ift\}}], \q x\in X.
\de
\end{lem}

\bg{proof}
We will use the second successive
approximation scheme of the minimal nonnegative solution \ct{c1, hzt}.
Let $$g^{(1)}(x)=P(x, A)+r(1)F^{1}(x, A),\q x\in X,$$
and inductively
$$g^{(n+1)}(x)=\int_{A^{c}}g^{(n)}(y)P(x, d y)+
r(n+1)F^{n+1}(x, A),\q n\geq 1.$$
Then we have
$$g^{(1)}(x)=\hat r(1)F^{1}(x, A).$$
Assume that
$g^{(n)}(x)=\hat r(n)F^{n}(x, A)$. Then
$$\aligned
g^{(n+1)}(x)&=\int_{A^{c}}\hat r(n)F^{n}(y, A)P(x, d y)+
r(n+1)F^{n+1}(x, A)\\
&=\hat r(n)F^{n+1}(x, A)+r(n+1)F^{n+1}(x, A)\\
&=\hat r(n+1)F^{n+1}(x, A).\\
\endaligned
$$
Hence
$$g^{*}(x)=\sum_{n=1}^{\ift}g^{(n)}(x)
=\sum_{n=1}^{\ift}\hat{r}(n)F^{n}(x, A)=
\E_x[\hat{r}(\tau_A)1_{\{\tau_A<\ift\}}]$$
is the minimal nonnegative solution of equation
(\ref{lsf}).
\end{proof}

\begin{cor}\lb{rmll}
$(1)$ For $A\in \mathcal{B}^{+}(X)$ and $\kp\geq1$,
\be\lb{sdf}
\E_x\l[\kp^{\tau_A}1_{\{\tau_A<\ift\}}\r]=
\kp\int_{A^{c}}\E_y\l[\kp^{\tau_A}1_{\{\tau_A<\ift\}}\r]P(x, d y)+\kp P(x, A).
\de
Moreover, $\l\{\E_x\l[\kp^{\sgm_A}1_{\{\sgm_A<\ift\}}\r], x\in X\r\}$
is the minimal nonnegative solution of the equations
\be\lb{bui}\aligned
\left\{\begin{array}{ll}
    g(x)=\kappa\int_{A^{c}} g(y)P(x, d y)+\kappa P(x, A) ,&  x\in A^{c};\\
    g(x)=1 ,&  x\in A.
\end{array}    \right.\\
\endaligned
\de

$(2)$ For $A\in \mathcal{B}^{+}(X)$ and
$\ell\in\N$,
\be\lb{sfd}\aligned
\E_x\l[(\tau_A+1)^{\ell}1_{\{\tau_A<\ift\}}\r]&=
\int_{A^{c}}\E_y\l[(\tau_A+1)^{\ell}1_{\{\tau_A<\ift\}}\r]P(x, d y)\\
&+P(x, A)
+\sum_{k=0}^{\ell-1}{\ell \choose k}
\E_x\l[\tau_A^{k}1_{\{\tau_A<\ift\}}\r].\\
\endaligned
\de
Moreover,
$\l\{\E_x\l[(\sgm_A+1)^{\ell}1_{\{\sgm_A<\ift\}}\r], x\in X\r\}$ is the
minimal nonnegative solution of the equations
\be\lb{bus1}\aligned
\left\{\begin{array}{ll}
    g(x)=\int_{A^{c}} g(y)P(x, d y)+P(x, A)+
    \sum_{k=0}^{\ell-1}{\ell \choose k}
    \E_x\l[\tau_A^{k}1_{\{\tau_A<\ift\}}\r],&  x\in A^{c};\\
    g(x)=1 ,&  x\in A.
\end{array} \right.\\
\endaligned
\de
\end{cor}
\bg{proof}
(1) Set $\hat{r}(n)=\kp^{n}$ with $\kp\geq1$ in Lemma \ref{r}. Then $$r(0)=\hat r(0)=1,\q r(n)=\kp^{n}-\kp^{n-1},\q n\geq1.$$
Hence
\be\lb{ssa}\aligned
&\E_x\l[\kp^{\tau_A}1_{\{\tau_A<\ift\}}\r]=
\int_{A^{c}}\E_y\l[\kp^{\tau_A}1_{\{\tau_A<\ift\}}\r]P(x, d y)\\&\q\q\q+P(x, A)
+\E_x\l[\l(\kp^{\tau_A}-\kp^{\tau_A-1}\r)1_{\{\tau_A<\ift\}}\r].\nnb\\
\endaligned
\de
Thus, \rf{sdf} holds by rearranging terms. Moreover,
by the localization theorem
and the comparison theorem of the minimal nonnegative solution
(see \ct[Chapter 2]{c1}),
the minimal nonnegative solution of equations \rf{bui} is $\l\{\E_x\l[\kp^{\sgm_A}1_{\{\sgm_A<\ift\}}\r], x\in X\r\}$.

(2) Set $\hat{r}(n)=(n+1)^{\ell}$ in Lemma \ref{r}. The proof is similar to that of (1), is then omitted.
\end{proof}

Now, we are ready to prove Theorem \ref{zzz}.

\bg{proof}[Proof of Theorem $\ref{zzz}$.]
We prove first $(3)\Rar(2)\Rar(1)\Rar(5)\Rar(3)$, and then $(2)\Leftrightarrow(4)$.

$(3)\Rar(2)$. Since $\sup_{x\in A}E_x\l[\kappa^{\tau_A}1_{\{\tau_A<\ift\}}\r]<\ift$,
we have for all $\dlt>0$,
there exists $N_0$ large enough such that
\be\lb{gis}
\sup_{x\in A}\sum_{n=N_0+1}^{\ift}\kappa^{n}F^{n}(x, A)\leq \dlt/2.
\de
Moreover, since $\sup_{x\in A}L(x, A)<1$,
there exists a constant $\lmd>1$ satisfying
\be\lb{gsi}
\sup_{x\in A}\sum_{n=1}^{N_0}\lmd^{n}F^{n}(x, A)\leq 1-\dlt.
\de
Set $\tld\kp=\min\{\kp, \lmd\}$. Then combining \rf{gis} with \rf{gsi}, we have
$$\sup_{x\in A}\sum_{n=1}^{\ift}\tld\kappa^{n}F^{n}(x, A)\leq
\sup_{x\in A}\sum_{n=1}^{N_0}\lmd^{n}F^{n}(x, A)+
\sup_{x\in A}\sum_{n=N_0+1}^{\ift}\kappa^{n}F^{n}(x, A)\leq
1-\dlt/2<1.$$

$(2)\Rar(1)$ and $(1)\Rar(5)$ follow from Lemmas \ref{n} and \ref{nc}(1), respectively.

$(5)\Rar(3)$.
Suppose that $\Phi$ is geometrically transient. Then
there exist $B\in\mathcal{B}^{+}(X)$ and $\kappa>1$
such that
  \be\lb{mm}
  \sup_{x\in B}\sum_{n=1}^{\ift}\kp^{n}P^{n}(x, B)<\ift.
  \de
Since $B\in\mathcal{B}^{+}(X)$, it follows from Corollary \ref{y61} and Proposition \ref{ti8} that there exists a petite set $A\subset B$
with $\psi(A)>0$ such that $$\sup_{x\in A}L(x, A)<1.$$
On the other hand, noting that $A\subset B$ and $F^{n}(x, A)\leq P^{n}(x, A)$ for all $x\in X$,
we get from \rf{mm} that
$$\sup_{x\in A}\sum_{n=1}^{\ift}\kp^{n}F^{n}(x, A)\leq\sup_{x\in A}\sum_{n=1}^{\ift}\kp^{n}P^{n}(x, A)\leq\sup_{x\in B}\sum_{n=1}^{\ift}\kp^{n}P^{n}(x, B)<\ift.$$

$(4)\Rar(2)$.
If \rf{34e} holds with $A=X$, then
$P(x, X)\leq P W(x)\leq b$ for $x\in X$,
hence for $1<\kp<b^{-1}$,
$$\sup_{x\in X}\E_x\l[\kp^{\tau_X}1_{\{\tau_X<\ift\}}\r]=\sup_{x\in X}\kp P(x, X)\leq\kp b<1.$$

Suppose that (\ref{34e}) holds with $A\not=X$
and $b<\lmd$. Then $W$ satisfies
\be\lb{dri7}\aligned
   \left\{\begin{array}{ll}
   W(x)\geq\lmd^{-1}P W(x)\geq\lmd^{-1}\int_{A^{c}}W(y)P(x, d y)+\lmd^{-1}P(x, A),
   &  x\in A^{c};\\
   W(x)\geq1,&  x\in A.
   \end{array}\right.\nnb\\
\endaligned
\de
According to $(\ref{bui})$, the minimal nonnegative solution of
the inequalities is given by $\E_x\l[\lmd^{-\sigma_A}1_{\{\sgm_A<\ift\}}\r]$,
hence
\begin{equation*}
\E_x\l[\lmd^{-\sigma_A}1_{\{\sgm_A<\ift\}}\r]\leq W(x),\q x\in A^{c}.
\end{equation*}
Combining this inequality with \rf{sdf}, and
noting that $P W(x)\leq b<\lmd$ for $x\in A$, we obtain that for $x\in A$,
\begin{equation*}\aligned
\E_x\l[\lmd^{-\tau_A}1_{\{\tau_A<\ift\}}\r]&=
\lmd^{-1}\int_{A^{c}}\E_y\l[\lmd^{-\sigma_A}1_{\{\sgm_A<\ift\}}\r]
P(x, d y)+\lmd^{-1}
P(x, A)\\
&\leq \lmd^{-1}\int_{A^{c}}W(y)P(x, d y)+\lmd^{-1}P(x, A)\\
&\leq\lmd^{-1}\l[-\int_{A}W(y)P(x, d y)+b\r]+\lmd^{-1}P(x, A)\\
&\leq\lmd^{-1}b<1.\\
\endaligned
\end{equation*}
Thus, \rf{xi01} holds with $\kappa=\lambda^{-1}$.

If $\lmd\leq b<1$, then there exists $\varepsilon>0$ such that
$\lmd<b+\varepsilon<1$, and $W$ satisfies
\be\lb{dri1}\aligned
   \left\{\begin{array}{ll}
     W(x)>(b+\varepsilon)^{-1}P W(x)\geq(b+\varepsilon)^{-1}\int_{A^{c}}
     W(y)P(x, d y)+(b+\varepsilon)^{-1}P(x, A),
     &  x\in A^{c};\\
     W(x)\geq1,&  x\in A.\nonumber
   \end{array}    \right.\\
\endaligned
\de
Using a similar argument, we have
\begin{equation*}
\E_x\l[(b+\varepsilon)^{-\sigma_A}1_{\{\sgm_A<\ift\}}\r]\leq W(x),\q x\in A^{c},
\end{equation*}
and for $x\in A$,
\begin{equation*}\aligned
\E_x\l[(b+\varepsilon)^{-\tau_A}1_{\{\tau_A<\ift\}}\r]&
=(b+\varepsilon)^{-1}\int_{A^{c}}
\E_y\l[(b+\varepsilon)^{-\sigma_A}1_{\{\sgm_A<\ift\}}\r]P(x, d y)\\
&+(b+\varepsilon)^{-1}P(x, A)\\
&\leq(b+\varepsilon)^{-1}\int_{A^{c}}W(y)P(x, d y)+
(b+\varepsilon)^{-1}P(x, A)\\
&\leq(b+\varepsilon)^{-1}\l[-\int_{A}W(y)P(x, d y)+
b\r]+(b+\varepsilon)^{-1}P(x, A)\\
&=(b+\varepsilon)^{-1}b<1.\\
\endaligned
\end{equation*}
Then \rf{xi01} holds with $\kappa=(b+\varepsilon)^{-1}$.

$(2)\Rar(4)$.
Set $W(x)=\mathbb{E}_x\l[\kp^{\sgm_{A}}1_{\{\sgm_A<\ift\}}\r]$ for $x\in X$.
Then by Corollary \ref{rmll}(1),
\begin{equation*}
P W(x)=\kappa^{-1}W(x),\q x\in A^{c},
\end{equation*}
and
\begin{equation*}
P W(x)=\kappa^{-1}\mathbb{E}_x\l[\kappa^{\tau_{A}}1_{\{\tau_A<\ift\}}\r]\leq
\kappa^{-1}\sup_{x\in A}\mathbb{E}_x\l[\kappa^{\tau_{A}}1_{\{\tau_A<\ift\}}\r],\q x\in A.
\end{equation*}
Thus, \rf{34e} holds with $\lmd=\kappa^{-1}$ and
$b=\kappa^{-1}\sup_{x\in A}\mathbb{E}_x\l[\kappa^{\tau_{A}}1_{\{\tau_A<\ift\}}\r]$.
\end{proof}

In the drift condition \rf{34e}, we require that
$P W(x)\leq b<1$ for all $x\in A$, which
is sometimes difficult to apply.
Hence we provide the following more practical drift condition.

(\textbf{GT})\q There exist some petite set $A\in\mathcal{B}^{+}(X)$, a constant $\lmd\in(0, 1)$, a nonnegative

\q\q\q\;\;function $W(x)<1$ for $x\in A^{c}$ and $W(x)\geq1$ for $x\in A$,
satisfying the drift

\q\q\q\;\;condition
$$P W(x)\leq\lmd W(x),\q x\in A^{c}.$$
In Remark \ref{s683}(2), we point out that
if the drift condition \rf{34e} holds with
a stochastic transition kernel $P$, then $\{x\in A^{c}: W(x)<1\}\not=\emptyset$.
Here, we strengthen the condition as $W(x)<1$ for all $x\in A^{c}$.

\bg{thm}
Suppose that $\Phi$ is a $\psi$-irreducible chain. If condition \mbox{(\textbf{GT})}
holds, then $\Phi$ is geometrically transient.
\end{thm}
\bg{proof}
Suppose that (\textbf{GT}) holds. Then we have
\be\lb{drw6}\aligned
   \left\{\begin{array}{ll}
     P W(x)\leq \lmd W(x),
     &  x\in A^{c};\\
     W(x)\geq1,&  x\in A.\nonumber
   \end{array}    \right.\\
\endaligned
\de
Thus, by Corollary \ref{rmll}(1), we get for $x\in A^{c}$,
\be\lb{cy65}
\E_x\l[\lmd^{-\tau_A}1_{\{\tau_A<\ift\}}\r]\leq W(x)<1.
\end{equation}
Hence
$$L(x, A)\leq\E_x\l[\lmd^{-\tau_A}1_{\{\tau_A<\ift\}}\r]<1,\q x\in A^{c},$$
and
$$L(x, A)=\int_{A^{c}}L(y, A)P(x, d y)+P(x, A)<1,\q x\in A.$$
Then there exist $\dlt<1$ and $B\subset A$ with $\psi(B)>0$ such that $L(x, B)\leq \dlt$
for all $x\in B$. That is,
\be\lb{by48}
\sup_{x\in B}L(x, B)<1.
\de

On the other hand, by \rf{cy65}, we have
\begin{equation*}\aligned
\E_x\l[\lmd^{-\tau_A}1_{\{\tau_A<\ift\}}\r]&
=\lmd^{-1}\int_{A^{c}}
\E_y\l[\lmd^{-\tau_A}1_{\{\tau_A<\ift\}}\r]P(x, d y)+\lmd^{-1}P(x, A)\\
&<\lmd^{-1}P(x, A^{c})+
\lmd^{-1}P(x, A)\leq\lmd^{-1},\q x\in A.\\
\endaligned
\end{equation*}
That is, \be\lb{sr54}
\sup_{x\in A}\E_x\l[\lmd^{-\tau_A}1_{\{\tau_A<\ift\}}\r]
=:b<\ift.
\end{equation}
In the following, we will prove that for some $r>1$,
$$\sup_{x\in B}\E_x\l[r^{\tau_B}1_{\{\tau_B<\ift\}}\r]<\ift.$$
This together with \rf{by48} yields the desired assertion.
The proof can be divided into three steps.

(a) First, we prove
\be\lb{ht5}
\sup_{x\in A}\E_x\l[\sum_{k=0}^{\tau_A-1}r^{k}\sum_{n=1}^{\ift}\lmd^{-n}F^{n}(\Phi_{k}, A)\r]<\ift,\q 1<r<\lmd^{-1}.
\end{equation}
Set
\be\lb{drx4}\aligned
  f(x)=\left\{\begin{array}{ll}
     \sum_{n=1}^{\ift}\lmd^{-n}F^{n}(x, A),
     &  x\in A^{c};\\
     1,&  x\in A.
   \end{array}    \right.\\
\endaligned
\de
Then $f$ satisfies
\begin{equation*}\aligned
P f(x)&=\lmd f(x)1_{A^{c}}(x)+\lmd\sum_{n=1}^{\ift}\lmd^{-n}F^{n}(x, A)1_{A}(x)\\
&\leq r^{-1}f(x)1_{A^{c}}(x)-\varepsilon f(x)1_{A^{c}}(x)+\lmd b 1_{A}(x),\\
\endaligned
\end{equation*}
for all $x\in X$ and $1<r<\lmd^{-1}$, where $\varepsilon=r^{-1}-\lmd$.
By defining $Z_k=r^{k}f(\Phi_k)$ for $k\in\Z_+$, it follows that
\begin{equation*}\aligned
&\E\l[Z_{k+1}|\mathcal{F}_k\r]=r^{k+1}\E\l[f(\Phi_{k+1})|\mathcal{F}_k\r]\\
&\leq r^{k+1}\l[r^{-1}f(\Phi_{k})1_{A^{c}}(\Phi_{k})-\varepsilon f(\Phi_{k})1_{A^{c}}(\Phi_{k})+\lmd b 1_{A}(\Phi_{k})\r]\\
&\leq Z_k-\varepsilon r^{k+1}f(\Phi_{k})1_{A^{c}}(\Phi_{k})+\lmd b r^{k+1}1_{A}(\Phi_{k}).\\
\endaligned
\end{equation*}
Then by \ct[Proposition 11.3.2]{meyn}, for all $C\in\mathcal{B}^{+}(X)$,
$$\E_x\l[\sum_{k=0}^{\tau_C-1}\varepsilon r^{k+1}f(\Phi_{k})1_{A^{c}}(\Phi_{k})\r]\leq
Z_0(x)+\E_x\l[\sum_{k=0}^{\tau_C-1}\lmd b r^{k+1}1_{A}(\Phi_{k})\r].$$
Multiplying by $\varepsilon^{-1}r^{-1}$ and noting that $Z_0(x)=f(x)$, we obtain that
$$\E_x\l[\sum_{k=0}^{\tau_C-1}r^{k}f(\Phi_{k})1_{A^{c}}(\Phi_{k})\r]\leq
\varepsilon^{-1}r^{-1}f(x)+\varepsilon^{-1}\lmd b\E_x\l[\sum_{k=0}^{\tau_C-1} r^{k}1_{A}(\Phi_{k})\r],$$
which yields that $$\E_x\l[\sum_{k=0}^{\tau_A-1}r^{k}f(\Phi_{k})1_{A^{c}}(\Phi_{k})\r]\leq
\varepsilon^{-1}r^{-1}f(x)+\varepsilon^{-1}\lmd b1_{A}(x).$$
Thus, by \rf{drx4},
$$\sup_{x\in A}\E_x\l[\sum_{k=1}^{\tau_A-1}r^{k}\sum_{n=1}^{\ift}\lmd^{-n}F^{n}(\Phi_k, A)\r]\leq
\varepsilon^{-1}r^{-1}+\varepsilon^{-1}\lmd b<\ift.$$
Combining this with \rf{sr54}, we get \rf{ht5}.

(b) Noting that $A$ is petite,
according to \rf{ht5} and the proof of \ct[Theorem 15.2.1]{meyn},
we obtain
$$\sup_{x\in A}\E_x\l[\sum_{k=0}^{\tau_C-1}r^{k}\sum_{n=1}^{\ift}\lmd^{-n}F^{n}(\Phi_{k}, A)\r]<\ift,\q C\in\mathcal{B}^{+}(X).$$

(c) For all $C\in\mathcal{B}^{+}(X)$, by the Markov property and noting that $\lmd<1$, we have
\begin{equation*}\aligned
&\E_x\l[\sum_{k=0}^{\tau_C-1}r^{k}\sum_{n=1}^{\ift}\lmd^{-n}F^{n}(\Phi_k, A)\r]
\geq \E_x\l[\sum_{k=0}^{\tau_C-1}r^{k}\sum_{n=1}^{\ift}F^{n}(\Phi_k, A)\r]\\
&=\E_x\l[\sum_{k=0}^{\tau_C-1}r^{k}\E_{\Phi_k}1_{\{\tau_A<\ift\}}\r]
=\E_x\l[\sum_{k=0}^{\ift}r^{k}1_{\{\tau_C\geq k+1\}}\E_{\Phi_k}1_{\{\tau_A<\ift\}}\r]\\
&=\E_x\l[\sum_{k=0}^{\ift}r^{k}1_{\{\tau_C\geq k+1\}}\E\l(1_{\{\theta^{k}\tau_A<\ift\}}|\mathcal{F}_k\r)\r]\\
&=\E_x\l[\sum_{k=0}^{\ift}r^{k}\E\l(1_{\{\tau_C\geq k+1\}}1_{\{\theta^{k}\tau_A<\ift\}}|\mathcal{F}_k\r)\r]\\
&=\E_x\l[\sum_{k=0}^{\ift}r^{k}1_{\{\tau_C\geq k+1\}}1_{\{\theta^{k}\tau_A<\ift\}}\r]\\
&=\E_x\l[\sum_{k=0}^{\tau_C-1}r^{k}1_{\{\theta^{k}\tau_A<\ift\}}\r],\\
\endaligned
\end{equation*}
where $\theta$ is the usual shift operator. It follows from (b) that
$$\sup_{x\in A}\E_x\l[\sum_{k=0}^{\tau_C-1}r^{k}1_{\{\theta^{k}\tau_A<\ift\}}\r]<\ift,\q C\in\mathcal{B}^{+}(X).$$
Noting that $B\subset A$, we arrive at
$$\sup_{x\in B}\E_x\l[\sum_{k=0}^{\tau_B-1}r^{k}1_{\{\theta^{k}\tau_B<\ift\}}\r]=
\sup_{x\in B}\E_x\l[\sum_{k=0}^{\tau_B-1}r^{k}1_{\{\tau_B<\ift\}}\r]<\ift,$$
which yields that
$$\sup_{x\in B}\E_x\l[r^{\tau_B}1_{\{\tau_B<\ift\}}\r]<\ift.$$
\end{proof}

\subsection{Strongly geometric transience}\label{tt}

In the previous section, we considered the geometric transience which
satisfies
$$\sup_{x\in A}\sum_{n=1}^{\ift}\kp^{n}P^{n}(x, A)<\ift,$$
for some $A\in\mathcal{B}^{+}(X)$ and $\kp>1$. However, in practice,
there exist a great deal of chains with sub-stochastic transition kernel,
for which we can further study the strongly geometric transience.
\begin{defn}\lb{uio}
The chain $\Phi$ is called strongly geometric transient
if there exists a constant $\kp>1$ such that
\be\lb{jit}
\sum_{n=1}^{\ift}\kp^{n}P^{n}(x, X)<\ift,\q x\in X.
\end{equation}
\end{defn}

Suppose that \rf{jit} holds and set
$A_i=\l\{x\in X: \sum_{n=1}^{\ift}\kp^{n}P^{n}(x, X)\leq i\r\}$ for $i\in\N$.
Then we have
$$X=\bigcup_{i=1}^{\ift}A_i\q\mbox{and}\q \sup_{x\in A_i}\sum_{n=1}^{\ift}\kp^{n}P^{n}(x, A_i)<\ift,\q i\geq1.$$
This implies that
if $\Phi$ is strongly geometric transient, it is geometrically transient.
For the converse, let $P=(p_{i j})$ be a transition kernel on
$X=\N$ with
\begin{equation}\lb{765}
P=\left(
\begin{array}{ccccc}
 0      & \gm_1         \\
 \bt_2  & 0      & \gm_2   \\
 \bt_3  & 0      & 0      & \gm_3   \\
 \vdots & \vdots & \vdots & \vdots  & \ddots \\
\end{array}
\right).
\end{equation}
If $\gm_1=1$, $\gm_k=(k-1)/k$ and $\bt_k=4^{-k}$
for $k\geq 2$, then the chain is geometrically
transient but not strongly geometric transient
by Theorems \ref{zzz}(2) and \ref{28}(2), respectively.

Since the transition kernel of strongly geometric transient chain
is sub-stochastic,
there is a positive
probability that the chain can ``escape" from $X$.
Let
\be\lb{876}
\tau=\sup\{n\geq0: \Phi_n\in X\}.
\de
Then for all bounded measurable function $f$ on $X$,
\be\lb{h}
P^{n}f(x)=\mathbb{E}_x
\left[f(\Phi_{n})1_{\{\tau>n\}}\right],\q x\in X.
\de

For the strongly geometric transience, we have the
following criteria.
\begin{thm}\lb{28}
Assume that $\Phi$ is $\psi$-irreducible, and
$\P_x\{\tau<\ift\}=1$ for all $x\in X$.
Then the following statements are equivalent.

$(1)$ The chain $\Phi$ is strongly geometric transient.

$(2)$ For $x\in X$, there exists a constant $\kp>1$ such that
$E_x \l[\kp^{\tau}\r]<\ift$.

$(3)$ There exist some constant $\lmd\in(0, 1)$
and a finite function $W\geq 1$ such that
\be\lb{zse}
P W(x)\leq \lmd W(x),\q x\in X.
\de
\end{thm}
\bg{proof}
$(1)\Rar(2)$. Suppose that \rf{jit} holds.
According to (\ref{h}), we have
\be\lb{j}
P^{n}(x, X)=\sup_{|f|\leq1}P^{n}f(x)=
\sup_{|f|\leq1}\mathbb{E}_x
\left[f(\Phi_{n})1_{\{\tau>n\}}\right]=
\mathbb{P}_x\{\tau>n\},\q x\in X.\nnb
\de
Hence for $x\in X$,
\begin{eqnarray*}
\begin{aligned}
&\mathbb{E}_x \l[\kp^{\tau}\r]=
(\kp-1)\sum_{m=0}^{\ift}\kp^{m}\mathbb{P}_x\{\tau>m\}+1\\
&=(\kp-1)\sum_{m=0}^{\ift}\kp^{m}P^{m}(x, X)+1<\ift.\\
\end{aligned}
\end{eqnarray*}

$(2)\Rar(3)$.
Denote by
$\widehat{X}=X\cup\{\partial\}$ the one point compactification of $X$. Let $\mathcal{B}(\hat X)=\sgm\l(\mathcal{B}(X)\cup\{\partial\}\r)$, and
\be\lb{89}\aligned\widehat{P}(x, A)=
    \left\{\begin{array}{ll}
     P(x, A),    & x\in X,\, A\in \mathcal{B}(X);\\
     1-P(x, X),  & x\in X,\, A\in\mathcal{B}(\hat X)\setminus\mathcal{B}(X);\\
     1,          & x\in \{\partial\},\, A\in\mathcal{B}(\hat X)
     \setminus\mathcal{B}(X);\\
     0,          & x\in \{\partial\},\, A\in \mathcal{B}(X).\\
   \end{array}    \right.\nnb\\
   \endaligned
   \de
Then $\hat P$ is a stochastic transition kernel. For $x\in \hat X$ and $\kp>1$, set
$$\hat{W}(x)=\E_x\l[\kp^{\tau}\r]1_{\l\{x\in X\r\}}+1_{\l\{x\in \{\partial\}\r\}}.$$
Then by Corollary \ref{rmll}(1), $\hat W$ satisfies
\be\lb{81s}\aligned
    \left\{\begin{array}{ll}
     \widehat{P} \widehat{W}(x)=\kp^{-1}\widehat{W}(x), & x\in X;\\
     \widehat{W}(\{\partial\})=1.\\
   \end{array}    \right.\nnb\\
   \endaligned
   \de
Let $W(x)=\widehat{W}(x)$ for $x\in X$. Then
$W\geq 1$ and
$$\aligned
&\kp^{-1} W(x)=\int_{X}\widehat{W}(y)\widehat{P}(x, d y)
+\int_{\{\partial\}}\widehat{W}(y)\widehat{P}(x, d y)\\
&=\int_{X}W(y)P(x, d y)+\widehat{P}(x, \{\partial\})\geq P W(x),\\
\endaligned
$$
which finishes the proof by letting $\lmd=\kp^{-1}$.

$(3)\Rar(1)$. Iterating the inequality (\ref{zse}) and noting that $W\geq1$,
we have
$$ P^{n}(x, X)\leq P^{n}W(x)\leq\lmd^{n}W(x),\q n\geq1.$$
Thus, $(1)$ holds with $1<\kp<\lmd^{-1}$.
\end{proof}

In the next, we will study the $V$-uniform transience for all function $V\geq1$,
which is closely related to the strongly geometric transience.

\begin{defn}\lb{yq}
The chain $\Phi$ is called V-uniformly transient for $V\geq1$, if
\be\lb{lkj}
||P^{n}||_{V}:=\sup_{x\in X}\frac{P^{n}V(x)}{V(x)}\rightarrow 0,
\q n\rar\ift.
\de
\end{defn}
Since $||\cdot||_{V}$ is an operator norm,
$||P^{m+n}||_{V}\leq||P^{m}||_{V}||P^{n}||_{V}$ for $m$, $n\in\Z_{+}$.
Thus, the convergence rate in \rf{lkj} must be geometric.

\begin{thm}\lb{21}
Assume that $\Phi$ is a $\psi$-irreducible chain.
Then the following statements are equivalent.

$(1)$ The chain $\Phi$ is $V$-uniformly transient for some $V\geq1$.

$(2)$ There exist some constant $\lmd\in(0, 1)$
and a finite function $W\geq 1$ such that
\be\lb{zdf}
P W(x)\leq \lmd W(x), \q x\in X,\nnb
\de
where $W$ is equivalent to $V$ in the sense that
$c^{-1}V\leq W\leq c V $
for some $c\geq 1$.
\end{thm}
\bg{proof}
$(1)\Rightarrow(2)$.
Suppose that there exist constants $R<\ift$ and $\rho<1$ such that
$||P^{n}||_{V}\leq R \rho^{n}$ for $n\geq0$.
Then since $\rho<1$, there exists $n_0\in\N$ large enough such that $R \rho^{n_0}<\beta^{-1}$ for some $\beta>1$. Set
$$W(x)=\sum_{i=0}^{n_0-1}\beta^{\frac{i}{n_0}}P^{i}V(x).$$
Then noting that
$P^{n}V(x)\leq R \rho^{n}V(x)$, we have
$$V(x)\leq W(x)\leq\sum_{i=0}^{n_0-1}
\beta^{\frac{i}{n_0}}R \rho^{i}V(x)\leq \beta n_0 R V(x),\q x\in X.$$
Moreover, in view of $R\rho^{n_0}<\beta^{-1}$, we get
$$\aligned
P W(x)&=\sum_{i=0}^{n_0-1}\beta^{\frac{i}{n_0}}P^{i+1}V(x)
=\sum_{i=1}^{n_0}\beta^{\frac{i-1}{n_0}}P^{i}V(x)\\
&=\beta^{-\frac{1}{n_0}}\sum_{i=1}^{n_0-1}\beta^{\frac{i}{n_0}}P^{i}V(x)
+\beta^{1-\frac{1}{n_0}}P^{n_0}V(x)\\
&\leq\beta^{-\frac{1}{n_0}}\sum_{i=1}^{n_0-1}\beta^{\frac{i}{n_0}}P^{i}V(x)
+\beta^{1-\frac{1}{n_0}}R\rho^{n_0}V(x)\\
&\leq\beta^{-\frac{1}{n_0}}\sum_{i=1}^{n_0-1}\beta^{\frac{i}{n_0}}P^{i}V(x)
+\beta^{-\frac{1}{n_0}}V(x)=\beta^{-\frac{1}{n_0}}W(x),\\
\endaligned
$$
which yields the conclusion by
letting $\lmd=\beta^{-\frac{1}{n_0}}$.

$(2)\Rightarrow(1)$. Since
$P^{n} W\leq \lmd^{n}W$ and
$c^{-1}V\leq W\leq c V$, we obtain that
$$||P^{n}||_{V}\leq
\sup_{x\in X}\frac{c P^{n}W(x)}{c^{-1}W(x)}\leq c^{2}\lmd^{n}\rightarrow 0,\q n\rightarrow\ift.$$
\end{proof}

\subsection{Uniformly geometric transience}\label{ut}

In this section, we will study a stronger type of geometric transience,
which the convergence in \rf{jit} is uniform with respect to initial states..

\begin{defn}\lb{uiio}
The chain $\Phi$ is called uniformly geometric transient
if there exists a constant $\kp>1$ such that
\bg{equation}\lb{cy7}
\sup_{x\in X}\sum_{n=1}^{\ift}\kp^{n}P^{n}(x, X)<\ift.
\end{equation}
\end{defn}

Obviously, uniformly geometric transient chains are
strongly geometric transient, but not conversely.
In fact, let $P=(p_{i j})$ be a random walk on $\Z_+$ with $p_{j, j+1}=p$,
$p_{j, j-1}=1-p$ and $p_{j j}=0$ for $j\geq0$. If $p<1/2$, then
it is strongly geometric transient but not uniformly
geometric transient according to Theorems \ref{28}(2) and \ref{39}(5), respectively.

Let $\tau$ be defined by \rf{876}. Then for the uniformly geometric transience, we obtain the following results.
\begin{thm}\lb{39}
Assume that $\Phi$ is $\psi$-irreducible, and
$\P_x\{\tau<\ift\}=1$ for all $x\in X$.
Then the following statements are equivalent.

$(1)$ The chain is uniformly geometric transient.

$(2)$ There exists a constant $\kp>1$ such that
$\sup_{x\in X}E_x \l[\kp^{\tau}\r]<\ift$.

$(3)$ There exists some constant $\lmd\in(0, 1)$ and
a bounded function $W\geq 1$ such that
$P W(x)\leq \lmd W(x)$ for $x\in X$.

$(4)$ There exists some $n_0\in\N$ such that $\sup_{x\in X}P^{n_0}(x, X)<1$.

$(5)$
$\sup_{x\in X}E_x\tau<\ift$.
\end{thm}
\bg{proof}
The proof of $(1)\Leftrightarrow(2)\Leftrightarrow(3)$ is similar to that
of Theorem \ref{28}. Thus, we only prove $(1)\Leftrightarrow(4)$
and $(2)\Leftrightarrow(5)$.

$(1)\Rar(4)$.
Set $M=\sup_{x\in X}\sum_{n=1}^{\ift}\kp^{n}P^{n}(x, X)$. Then we have
$$\sup_{x\in X}P^{n}(x, X)\leq\kp^{-n}M,\q n\in\N.$$
Since $\kp>1$, there exists $n_0$ large enough
such that $\kp^{-n_0}M<1$, which implies (4) holds.

$(4)\Rar(1)$. Set $\dlt=\sup_{x\in X}P^{n_0}(x, X)$.
Then it is easy to obtain that
\be\lb{iu8}
\sup_{x\in X}P^{k n_0}(x, X)\leq\dlt^{k}, \q k\in\N.
\de
For $n\in\N$, write $n=k n_0+s$, where $k$ is the integer
part of $n/n_0$ and $0\leq s< n_0$. Then by \rf{iu8},
\begin{eqnarray*}
\begin{aligned}
&P^{n}(x, X)=\int_{X}P^{k n_0}(y, X)P^{s}(x, d y)\\
&\leq\sup_{y\in X}P^{k n_0}(y, X)P^{s}(x, X)\leq \dlt^{k}\leq\dlt^{\fr{n-n_0}{n_0}}.\\
\end{aligned}
\end{eqnarray*}
Thus, \rf{cy7} holds with $1<\kp<\dlt^{-1/n_0}$. This proves (1).

$(2)\Leftrightarrow(5)$. Set $M=\sup_{x\in X}E_x\tau$.
Then by the minimal nonnegative solution theory (cf. \ct[Theorem 6.3.4]{hzt}), we have $$\sup_{x\in X}E_x \l[\tau^{n}\r]\leq n!M^{n}.$$
Hence for $1<\kp<e^{1/M}$,
\begin{eqnarray*}
\begin{aligned}
\log\kp~\E_x\tau&\leq\mathbb{E}_x \l[\kp^{\tau}\r]=\E_x \l[e^{\tau\log\kp}\r]\\&=\sum_{n=0}^{\ift}
(n!)^{-1}(\log\kp)^{n}E_x\l[ \tau^{n}\r]\\
&\leq\sum_{n=0}^{\ift}(\log\kp)^{n}M^{n}=(1-M\log\kp)^{-1},\\
\end{aligned}
\end{eqnarray*}
which completes the proof.
\end{proof}

\section{Algebraic transience}\label{als}

In this section, we will study algebraic transience. First, let us begin with the definition.

\begin{defn}\lb{fhg}
For an integer $\ell\geq1$, a set $A\in \mathcal{B}^{+}(X)$ is called uniformly $\ell$-transient if
$$\sup_{x\in A}\sum_{n=1}^{\infty}n^{\ell}P^{n}(x, A)<\ift.$$
The chain $\Phi$ is called $\ell$-transient
if it is $\psi$-irreducible and $X$ can be covered $\psi$-a.e.
by a countable number of uniformly $\ell$-transient sets.
\end{defn}

Similarly, if $\Phi$ is algebraically transient, then it is transient.
Since $\lim_{n\rightarrow\ift}\kp^{n}/n^{\ell}=\ift$
for $\kp>1$ and $\ell\geq 1$,
geometrically transient chains are algebraically transient,
but not conversely. Let $P=(p_{i j})$ be defined by \rf{765}
with
$\gm_1=1$, $\gm_k=(k-1)/k$ and $\bt_k=(k-1)/k^{\zeta}$
for $k\geq 2$ and some integer $\zeta\geq3$. Then the chain is algebraically transient, but it is not geometrically transient.

For the algebraic transience, we have the following criteria connecting
the ``local" algebraic transience, the first return, the
drift condition and the algebraic transience.
\bg{thm}\lb{aaq}
Let $\ell\geq1$ be an integer. Suppose that the chain $\Phi$ is $\psi$-irreducible.
Then the following statements are equivalent.

$(1)$ There exists a set $A\in\mathcal{B}^{+}(X)$ such that
$$\sup_{x\in A}\sum_{n=1}^{\infty}n^{\ell}P^{n}(x, A)<\ift.$$

$(2)$ There exists a set $A\in\mathcal{B}^{+}(X)$ such that
\be\lb{cy4s}
\sup_{x\in A}L(x, A)<1,\q \sup_{x\in A}\E_x\l[\tau_A^{\ell}1_{\{\tau_A<\ift\}}\r]<\ift.
\end{equation}

$(3)$  There exist some set
$A\in \mathcal{B}^{+}(X)$, constants $d\in(0, \ift)$,
$b\in(0, 1)$, and nonnegative functions
$W_i$, $i=0, 1, \cdots, \ell$, with $W_i(x_0)<\ift$ for some $x_0\in X$, satisfying for $i=0, 1, \cdots, \ell$,
\be\lb{s5e2}\aligned
\left\{\begin{array}{ll}
P W_{i}(x)\leq W_{i}(x)-(\ell-i)W_{i+1}(x), & x\in A^{c};\\
W_{i}(x)\geq1 ,& x\in A;\\
P W_0(x)\leq d,& x\in A;\\
P W_{\ell}(x)\leq b, & x\in A,\\
\end{array} \right.\\
\endaligned
\de
where $W_{\ell+1}=0$.

$(4)$ The chain $\Phi$ is $\ell$-transient.
\end{thm}

\bg{rem}
$(1)$ The set $A\in\mathcal{B}^{+}(X)$ is a petite set.

$(2)$ Since $P W_{\ell}(x)\leq b$ holds with $b\in(0, 1)$ in \rf{s5e2},
the set $\{x\in A^{c}: W_{\ell}(x)<1\}\not=\emptyset$ when $P$ is stochastic.
\end{rem}

To prove the theorem, we need
the following lemma. It gives the condition on the first return time which ensures that a set is uniformly $\ell$-transient.

\begin{lem}\lb{ng00}
Let $A\in \mathcal{B}^{+}(X)$ and $\ell\in\N$. Suppose that
$$\sup_{x\in A}L(x, A)<1,\q
\sup_{x\in A}\sum_{n=1}^{\ift}n^{\ell}F^{n}(x, A)<\ift.$$
Then $$\sup_{x\in A}\sum_{n=1}^{\ift}n^{\ell}P^{n}(x, A)<\ift.$$
\end{lem}
\bg{proof}
Set
$\delta=\sup_{x\in A}L(x, A)$ and
$M=\sup_{x\in A}\sum_{n=1}^{\ift}n^{\ell}F^{n}(x, A)$.
Then for fixed $N\in\mathbb{N}$, it follows from (\ref{le}) and
the binomial theorem that for $x\in A$,
\be\lb{plm}\aligned
&\sum_{n=1}^{N}n^{\ell}P^{n}(x, A)=\sum_{n=1}^{N}n^{\ell}
F^{n}(x, A)+\sum_{n=1}^{N}\sum_{m=1}^{n-1}
\int_{A}P^{m}(x, d y)F^{n-m}(y, A)n^{\ell}\\
&=\sum_{n=1}^{N}n^{\ell}
F^{n}(x, A)+\int_{A}\sum_{m=1}^{N-1}
P^{m}(x, d y)\sum_{n=m+1}^{N}F^{n-m}(y, A)(m+n-m)^{\ell}\\
&=\sum_{n=1}^{N}n^{\ell}F^{n}(x, A)+
\int_{A}\sum_{m=1}^{N-1}m^{\ell}P^{m}(x, d y)
\sum_{n=1}^{N-m}F^{n}(y, A)\\
&+\int_{A}\sum_{m=1}^{N-1}P^{m}(x,d y)\sum_{n=1}^{N-m}n^{\ell}
F^{n}(y, A)+\sum_{k=1}^{\ell-1}{\ell \choose k}\int_{A}\sum_{m=1}^{N-1}m^{k}P^{m}
(x, d y)\sum_{n=1}^{N-m}n^{\ell-k}F^{n}(y, A)\\
&\leq M+\dlt\sum_{m=1}^{N}m^{\ell}P^{m}(x, A)+M\sum_{m=1}^{N}P^{m}(x, A)+
M \sum_{k=1}^{\ell-1}{\ell \choose k}\sum_{m=1}^{N}m^{k}P^{m}(x, A).\nnb\\
\endaligned
\de
That is, for $x\in A$,
\be\lb{pln}\aligned
\sum_{n=1}^{N}n^{\ell}P^{n}(x, A)&
\leq \frac{M}{1-\delta}\l[1+\sum_{n=1}^{N}P^{n}(x, A)
+\sum_{k=1}^{\ell-1}{\ell \choose k}\sum_{n=1}^{N}n^{k}P^{n}(x, A)\r].
\endaligned
\de
On the other hand, by Lemma \ref{n}, we have
$$\sum_{n=1}^{\ift}P^{n}(x, A)\leq\fr{\dlt}{1-\dlt},\q x\in A.$$
Then we get
$$\aligned
\sum_{n=1}^{\ift}n P^{n}(x, A)
\leq \frac{M}{1-\delta}\l[1+\sum_{n=1}^{\ift}P^{n}(x, A)\r]
\leq \fr{M}{(1-\dlt)^{2}},\q x\in A.
\endaligned
$$
Combining these two inequalities with \rf{pln}, and by the induction argument,
we complete the proof.
\end{proof}

Now, it is ready to prove Theorem \ref{aaq}.

\bg{proof}[Proof of Theorem \mbox{\ref{aaq}}.]
$(2)\Rar(1)$ and $(1)\Rar(4)$ follow from Lemmas \ref{ng00} and \ref{nc}(1), respectively. $(4)\Rar(2)$ is similar to $(5)\Rar(3)$ of Theorem \ref{zzz}. Thus, we need only prove $(2)\Leftrightarrow(3)$.

$(3)\Rar(2)$. If \rf{s5e2} holds with $A=X$,
then the proof is similar to that of Theorem \ref{zzz}.

Suppose that (\ref{s5e2}) holds with $A\not=X$. Then $W_{\ell}$ satisfies
\be\lb{plfr}\aligned
   \left\{\begin{array}{ll}
   W_{\ell}(x)\geq P W_{\ell}(x)\geq \int_{A^{c}}W_{\ell}(y)P(x, d y)+P(x, A),&  x\in A^{c};\\
   W_{\ell}(x)\geq1,&  x\in A.
   \end{array}\right.\nnb\\
\endaligned
\de
According to \rf{bui}, the minimal nonnegative solution of the inequalities is given by
$L(x, A)1_{A^{c}}(x)+1_{A}(x)$. Hence $L(x, A)\leq W_{\ell}(x)$ for $x\in A^{c}$,
and for $x\in A$,
\be\lb{09aa}\aligned
&L(x, A)=\int_{A^{c}}L(y, A)P(x, d y)+P(x, A)\\
&\leq\int_{A^{c}}W_{\ell}(y)P(x, d y)+P(x, A)\\
&\leq-\int_{A}W_{\ell}(y)P(x, d y)+b+P(x, A)\leq b.\\
\endaligned
\de
Since $W_{\ell-1}(x)\geq1$ for $x\in A$, we have
\be\lb{09}\aligned
    W_{\ell-1}(x)\geq P W_{\ell-1}(x)+W_{\ell}(x)\geq \int_{A^{c}}W_{\ell-1}(y)P(x, d y)
    +P(x, A)+W_{\ell}(x),\q  x\in A^{c}.\nnb
\endaligned
\de
Set $\ell=1$ in $(\ref{bus1})$. Then it yields that
$$\sum_{n=1}^{\ift}(n+1)F^{n}(x, A)=\int_{A^{c}}\sum_{n=1}^{\ift}(n+1)F^{n}(y, A)
P(x, d y)+P(x, A)+L(x, A),\q x\in A^{c}.$$
Noting that $L(x, A)\leq W_{\ell}(x)$ for $x\in A^{c}$, it follows from the comparison theorem that,
\be\lb{jye}
\sum_{n=1}^{\ift}(n+1)F^{n}(x, A)\leq W_{\ell-1}(x),\q x\in A^{c}.\nnb
\de
Suppose that for all $i\leq \ell-1$,
$$\sum_{n=1}^{\ift}(n+1)^{i}F^{n}(x, A)
\leq W_{\ell-i}(x),\q x\in A^{c}.$$
Then
\be\lb{ci3}\aligned
&\sum_{k=0}^{\ell-1}{\ell \choose k}\E_x\l[\tau_A^{k}1_{\{\tau_A<\ift\}}\r]=
\sum_{n=1}^{\ift}\sum_{k=0}^{\ell-1}{\ell \choose k}
n^{k}F^{n}(x, A)\\
&\leq\sum_{n=1}^{\ift}\ell(n+1)^{\ell-1}F^{n}(x, A)\leq\ell W_{1}(x),\q x\in A^{c}.\\
\endaligned
\de
Since $W_{0}(x)\geq1$ for $x\in A$, we have
\be\lb{02}\aligned
    W_{0}(x)\geq P W_{0}(x)+\ell W_{1}(x)\geq
    \int_{A^{c}}W_{0}(y)P(x, d y)+P(x, A)+\ell W_{1}(x),
    \q  x\in A^{c}.\nnb
\endaligned
\de
Thus, combining this with \rf{bus1} and \rf{ci3}, we have
$$\sum_{n=1}^{\ift}(n+1)^{\ell}F^{n}(x, A)
\leq W_{0}(x),\q x\in A^{c}.$$
Noting that $P W_{0}(x)\leq d$, it follows from \rf{sfd} that for $x\in A$,
\be\lb{cu69}\aligned
\sum_{n=1}^{\ift}n^{\ell}F^{n}(x, A)&=\int_{A^{c}}
\sum_{n=1}^{\ift}(n+1)^{\ell}F^{n}(y, A)P(x, d y)+P(x, A)\\
&\leq \int_{A^{c}}W_{0}(y)P(x, d y)+P(x, A)\\
&\leq -\int_{A}W_{0}(y)P(x, d y)+d+P(x, A)\leq d.\\
\endaligned
\de
Thus, (2) holds by \rf{09aa} and \rf{cu69}.

$(2)\Rar(3)$. Suppose that \rf{cy4s} holds. For $i=0, 1, \cdots, \ell$,
set $$W_i(x)=(\ell-i)!\,\E_x\l[(\sgm_A+1)^{\ell-i}1_{\{\sgm_A<\ift\}}\r],
\q x\in X.$$
Then by Corollary \ref{rmll}(2) and noting that
$$\sum_{k=0}^{i-1}
{i \choose k}n^{k}=(n+1)^{i}-n^{i}\geq(n+1)^{i-1},$$
we obtain that
$$\aligned
W_{i}(x)&=\int_{A^{c}}W_{i}(y)P(x, d y)
+(\ell-i)! P(x, A)\\&+(\ell-i)! \sum_{n=1}^{\ift}\sum_{k=0}^{\ell-i-1}
{\ell-i \choose k}n^{k}F^{n}(x, A)\\
&\geq P W_i(x)+(\ell-i) (\ell-i-1)!
\sum_{n=1}^{\ift}(n+1)^{\ell-i-1}F^{n}(x, A)\\
&=P W_i(x)+(\ell-i) W_{i+1}(x),\q x\in A^{c},\\
\endaligned
$$
$$\aligned
P W_{0}(x)&=\ell ! \int_{A^{c}}\sum_{n=1}^{\ift}(n+1)^{\ell}
F^{n}(y, A)P(x, d y)+\ell ! P(x, A)\\
&=\ell ! \sum_{n=1}^{\ift}(n+1)^{\ell}
F^{n+1}(x, A)+\ell ! P(x, A)\\
&\leq \ell ! \sup_{x\in A}\sum_{n=1}^{\ift}n^{\ell}
F^{n}(x, A)<\ift,\q x\in A,\\
\endaligned
$$
and
$$P W_{\ell}(x)=L(x, A)\leq\sup_{x\in A}L(x, A)<1,\q x\in A.$$
\end{proof}
The drift conditions in Theorem \ref{aaq} can be generated
from two drift conditions by using the following lemma, which is
just a modification of \ct[Lemma 3.5]{jr}.
\bg{lem}\lb{xi9q}
Suppose that there exist $A\in\mathcal{B}^{+}(X)$, constants $d\in[0, \ift)$, $\alpha\in(0, 1)$, and a function $V\geq1$ satisfying
$$P V(x)\leq V(x)-d V^{\alpha}(x),\q x\in A^{c}.$$
Then for every $0<\eta\leq 1$,
$$P V^{\eta}(x)\leq V^{\eta}(x)-\eta d V^{\alpha+\eta-1}(x),\q x\in A^{c}.$$
\end{lem}
\bg{thm}\lb{098}
Let $\ell\geq1$ be an integer. Suppose that $\Phi$ is $\psi$-irreducible.
Then $\Phi$ is $\ell$-transient if there exist
some set $A\in\mathcal{B}^{+}(X)$, constants
$d\in(0, \ift)$, $b\in(0, 1)$, non-negative functions $V(x)\geq1$ for $x\in X$,
$W(x)\geq 1$ for $x\in A$ and $W(x)\leq 1$ for $x\in A^{c}$,
satisfying for $x\in X$,
\be\lb{pl19}\aligned
   \left\{\begin{array}{ll}
   P V(x)\leq \l(V(x)-\ell V^{1-\fr{1}{\ell}}(x)\r)1_{A^{c}}(x)+d 1_{A}(x);\\
   P W(x)\leq W(x) 1_{A^{c}}(x)+b 1_{A}(x).
   \end{array}\right.\\
\endaligned
\de
\end{thm}
\bg{proof}
By Lemma \ref{xi9q}, we have
$$P V^{1-\fr{i}{\ell}}(x)\leq V^{1-\fr{i}{\ell}}(x)-(\ell-i)V^{1-\fr{i+1}{\ell}}(x),\q x\in A^{c},\; i=0, 1, \cdots, \ell-1.$$
Set $W_i=V^{1-\fr{i}{\ell}}$ for $i=0, 1, \cdots, \ell-1$, and $W_{\ell}=W$.
Then for $x\in A^{c}$ and $i=0, 1, \cdots, \ell-2$,
$$P W_{i}(x)\leq W_{i}(x)-(\ell-i)W_{i+1}(x),$$
$$P W_{\ell-1}(x)\leq W_{\ell-1}(x)-1\leq W_{\ell-1}(x)-W_{\ell}(x),$$
and
$$P W_{\ell}(x)\leq W_{\ell}(x).$$
For $x\in A$, we get $P W_{0}(x)\leq d<\ift$ and $P W_{\ell}(x)\leq b<1$.
Thus, the drift conditions \rf{s5e2} hold, which imply the $\ell$-transience.
\end{proof}

\section{Applications}\label{exan}

This section is devoted to applying our results
to the random walk on $\R_+$
and the skip-free chain on $\Z_+$.

\subsection{The random walk on the half line}\label{ra}

In this section, we illustrate the applicability of the drift conditions
\rf{34e} and \rf{pl19} by the random walk on the half line.

\bg{exm}\lb{x8a}
(The random walk on the half line).
Suppose that $\Phi=\{\Phi_n: n\in\Z_+\}$ is defined by choosing an arbitrary
distribution for $\Phi_0$ and taking
$$\Phi_{n+1}=\l(\Phi_n+U_{n+1}\r)^{+},\q n\in\Z_+,$$
where $\l(\Phi_n+U_{n+1}\r)^{+}=\max\l(0, \Phi_n+U_{n+1}\r)$, and $(U_n)$ is a sequence of i.i.d. random variables taking values in $\R$ with
$$\Gm(-\ift, x]=\P\l(U\leq x\r),\q x\in\R.$$
\end{exm}

We write RWHL for short. The RWHL has a wide range of application,
it is both a model for storage systems and a model for queueing systems.
For more details, one can refer to \ct[Chapter 3]{meyn}.

If $\Gm(-\ift, 0)>0$, then the RWHL is
$\dlt_0$-irreducible and all compact sets are petite. By considering
the motion of the chain after it reaches $\{0\}$, we see that it
is also $\psi$-irreducible, where
$$\psi(A)=\sum_n P^{n}(0, A)2^{-n},\q A\in\mathcal{B}(X).$$
According to \ct[Proposition 4.2.2]{meyn}, $\psi$ is the maximal irreducibility measure.

Set $\beta=\int_{-\ift}^{\ift}x\Gm(d x)$. Proposition 9.5.1 in \ct{meyn} shows that if $\beta>0$, then the RWHL is transient.
For the geometric ergodicity and algebraic ergodicity of the chain, see \ct{dfm, jr, meyn, tt}.
Here, we will first study the geometric transience
of the RWHL by using the drift condition \rf{34e}.

\bg{thm}\lb{bte7}
Assume that $\beta>0$ and there exists $\gm>0$ such that
$$\int_{-\ift}^{0}\exp(-\gm x)\Gm(d x)<\ift.$$
Then the RWHL is geometrically transient.
\end{thm}
\bg{proof}
Set $A=\{0\}\in\mathcal{B}^{+}(X)$ and choose $W(x)=\exp(-s x)$, where $s\in(0, \gm]$ is to be specified later.
Then $W(0)=1$ and
\be\lb{p9d1}\aligned
P W(0)&=\int W(y)P(0, d y)=\int \exp(-s y)P(0, d y)\\
&=\int_{0}^{\ift}\exp(-s y)\Gm(d y)+\Gm(-\ift, 0]\\
&=\int_{-\ift}^{\ift}\exp(-s y)\Gm(d y)-\int_{-\ift}^{0}\l(\exp(-s y)-1\r)\Gm(d y)\\
&\leq\int_{-\ift}^{\ift}\l(\exp(-s y)-1\r)\Gm(d y)+1.\nnb\\
\endaligned
\de
Since
$s\mapsto s^{-1}\l(\exp(-s y)-1\r)$ is decreasing, we get
$$s^{-1}\int_{-\ift}^{\ift}\l(\exp(-s y)-1\r)\Gm(d y)\rar -\beta<0,
\q s\downarrow 0.$$
Thus, we can choose and fix $s_0$ sufficiently small such that
$$\int_{-\ift}^{\ift}\l(\exp(-s_0 y)-1\r)\Gm(d y)=\xi<0.$$
Hence we have
\be\lb{h7a}
P W(0)\leq\xi+1<1.
\end{equation}
Now set $W(x)=\exp(-s_0 x)$. Then for $x\in A^{c}$,
\be\lb{i9d0}\aligned
\fr{P W(x)-W(x)}{W(x)}&=\int\l(\fr{W(y)}{W(x)}-1\r)P(x, d y)=\int\l(\exp(-s_0 y+s_0 x)-1\r)P(x, d y)\\
&=\int_{-x}^{\ift}\l(\exp(-s_0 y)-1\r)\Gm(d y)+\l(\exp(s_0 x)-1\r)\Gm(-\ift, -x]\\
&=\int_{-\ift}^{\ift}\l(\exp(-s_0 y)-1\r)\Gm(d y)+\int_{-\ift}^{-x}\l(\exp(s_0 x)-\exp(-s_0 y)\r)\Gm(d y)\\
&\leq\int_{-\ift}^{\ift}\l(\exp(-s_0 y)-1\r)\Gm(d y)=\xi<0.\nnb\\
\endaligned
\de
Thus, $P W(x)\leq(\xi+1)W(x)$ for all $x\in A^{c}$. Combining this with \rf{h7a},
the RWHL is geometrically transient according to the drift condition \rf{34e}.
\end{proof}

Next, we will investigate the algebraic transience of the RWHL
by using the drift condition \rf{pl19}.
\bg{thm}\lb{xu87}
Assume that $\beta>0$ and there exists some integer $\ell\geq2$ such that
$$\int_{-\ift}^{0}(-x)^{\ell}\Gm(d x)<\ift.$$
Then the RWHL is $\ell$-transient.
\end{thm}
\bg{proof}
Since $\beta>0$, there exists $x_0>0$ such that
$\int_{-x_0}^{\ift}y\Gm(d y)=:\xi>0$. Set $$V(x)=(c x+a)^{-\ell}+1,$$
where $a\in(0, \ift)$ and $c\in(0, 1]$ are to be specified later. Then for $x>x_0$,
\be\lb{i9t7}\aligned
&P V(x)-V(x)=\int\l(V(y)-V(x)\r)P(x, d y)\\
&=\int\l((c y+a)^{-\ell}-(c x+a)^{-\ell}\r)P(x, d y)\\
&=\int_{-x}^{\ift}\l((c y+c x+a)^{-\ell}-(c x+a)^{-\ell}\r)\Gm(d y)+\l(a^{-\ell}-(c x+a)^{-\ell}\r)\Gm(-\ift, -x]\\
&=\int_{-x_0}^{\ift}\l((c y+c x+a)^{-\ell}-(c x+a)^{-\ell}\r)\Gm(d y)\\&+
\int_{-x}^{-x_0}\l((c y+c x+a)^{-\ell}-(c x+a)^{-\ell}\r)\Gm(d y)+\l(a^{-\ell}-(c x+a)^{-\ell}\r)\Gm(-\ift, -x]\\
&\leq\int_{-x_0}^{\ift}\l((c y+c x+a)^{-\ell}-(c x+a)^{-\ell}\r)\Gm(d y)+a^{-\ell}\\
&\leq V^{1-\fr{1}{\ell}}(x)\l(\int_{-x_0}^{\ift}\fr{(c x+a)^{\ell}-(c y+c x+a)^{\ell}}{(c x+a)(c y+c x+a)^{\ell}\l(1+(c x+a)^{\ell}\r)^{1-\fr{1}{\ell}}}\Gm(d y)+a^{-\ell}\r).\\
\endaligned
\de
Let $$f(c, x, y)=\fr{(c x+a)^{\ell}-(c y+c x+a)^{\ell}}{(c x+a)(c y+c x+a)^{\ell}\l(1+(c x+a)^{\ell}\r)^{1-\fr{1}{\ell}}},\q x>x_0,\; y>-x_0.$$
In the following, we turn to bound $f(c, x, y)$.
Note that, for $x+y\geq0$ and $-y\geq0$,
$$a^{\ell-2}(c x+a)^{\ell-2}\leq (c y+c x+a)^{\ell-2}(-c y+a)^{\ell-2},$$
since
\be\lb{i9zi}\aligned
&\log(c x+a)-\log(c y+c x+a)=\int_{(c y+c x+a)}^{(c x+a)}\fr{1}{z}d z\\
&\leq\int_{a}^{-c y+a}\fr{1}{z}d z
=\log(-c y+a)-\log a.\nnb\\
\endaligned
\de
For $y\in(-x_0, 0)$, we then get
\be\lb{xizi}\aligned
&\l|(c x+a)^{\ell}-(c y+c x+a)^{\ell}\r|=
(c y+ c x+a-c y)^{\ell}-(c y+c x+a)^{\ell}\\
&\leq -\ell c y(c y+c x+a)^{\ell-1}
+\fr{\ell(\ell-1)}{2}c^{2}y^{2}(c y+ c x+a-c y)^{\ell-2}\\
&\leq -\ell c y(c y+c x+a)^{\ell-1}
+\fr{\ell(\ell-1)}{2}c^{2}y^{2}(c y+ c x+a)^{\ell-2}\l(\fr{-c y+a}{a}\r)^{\ell-2}.\nnb\\
\endaligned
\de
For $y\geq 0$, we obtain
\be\lb{xia0}\aligned
&\l|(c x+a)^{\ell}-(c y+c x+a)^{\ell}\r|=
-(c y+ c x+a-c y)^{\ell}+(c y+c x+a)^{\ell}\\
&\leq \ell c y(c y+c x+a)^{\ell-1}
-\fr{\ell(\ell-1)}{2}c^{2}y^{2}(c y+ c x+a)^{\ell-2}\leq
\ell c y(c y+c x+a)^{\ell-1}.\nnb\\
\endaligned
\de
Thus, collecting the above estimates, we arrive at
\be\lb{aida0}\aligned
\fr{\l|f(c, x, y)\r|}{c}&\leq\l(\fr{-\ell y}{a^{2}}+\fr{\ell(\ell-1)(-y+a)^{\ell}}{2 a^{\ell+1}}\r)1_{\{y\in(-x_0, 0)\}}+\fr{\ell y}{a^{2}}1_{\{y\in[0, \ift)\}}=: g(y).\nnb\\
\endaligned
\de
From this, and noting that $\int_{-\ift}^{0}(-y)^{\ell}\Gm(d y)<\ift$, we have
$$\int_{-x_0}^{\ift}g(y)\Gm(d y)<\ift.$$
Then by the dominated convergence theorem, we get
\be\lb{i9t8}
\fr{1}{c}\int_{-x_0}^{\ift}\fr{(c x+a)^{\ell}-(c y+c x+a)^{\ell}}{(c x+a)(c y+c x+a)^{\ell}\l(1+(c x+a)^{\ell}\r)^{1-\fr{1}{\ell}}}\Gm(d y)\rar -\fr{\ell\xi}{a^{2}(1+a^{\ell})^{1-\fr{1}{\ell}}},\q c\downarrow 0.
\end{equation}
Therefore, combining \rf{i9t7} with \rf{i9t8}, we can choose and fix $c_0$ small enough such that for $x>x_0$,
\be\lb{aiia0}\aligned
P V(x)-V(x)&\leq -\ell V^{1-\fr{1}{\ell}}(x)\l(\fr{\xi c_0}{a^{2}(1+a^{\ell})^{1-\fr{1}{\ell}}}-\fr{1}{\ell a^{\ell}}\r).\nnb\\
\endaligned
\de
Hence we can also choose and fix $a_0$ sufficiently small such that
$$P V(x)-V(x)\leq-\ell V^{1-\fr{1}{\ell}}(x),\q x>x_0.$$
Since $V(x)$ is bounded for all $x\geq0$, it is obvious that $P V(x)\leq d<\ift$
for $x\leq x_0$.

Set $A=[0, x_0]$ and
$$W(x)=\fr{s x_0+1}{s x+1}1_{A^{c}}(x)+1_A(x).$$
Similarly, we can specify the constant $s>0$, and prove $P W(x)\leq W(x) 1_{A^{c}}(x)+b 1_{A}(x)$ with $b\in(0, 1)$. Then the chain is $\ell$-transient by Theorem \ref{098}.
\end{proof}

\subsection{The skip-free chain on $\Z_+$}\label{exa}

In this section, we will study the skip-free chain, and give the explicit criteria
for geometric transience and algebraic transience. For the ergodicity of the skip-free chain, see \ct{c1, m2, zyh1, zyh} and references within.

The remainder of the section is organized as follows.
In Theorem \ref{miw}, geometric transience is studied by
using the drift condition \rf{34e}. Through the first return time criteria \rf{cy4s}, the algebraic transience is discussed in Theorem \ref{mzj1}.
Finally, in Theorem \ref{bwo}, we consider the strongly geometric transience and uniformly geometric transience for the sub-stochastic skip-free chain.

\begin{exm}\lb{zzc}
(The skip-free Markov chain).
Let $P=(p_{i j})$ be an irreducible stochastic transition kernel on $X=\mathbb{Z}_{+}$ with
$p_{i j}=0$ for $j-i\geq2$.
\end{exm}

For $0\leq i<n$, define $p_{n}^{(i)}=\sum_{k=0}^{i}p_{n k}$, and inductively
\be\lb{000}
F_{n}^{(n)}=1, \q F_{n}^{(i)}=
\sum_{k=i}^{n-1}\fr{p_{n}^{(k)}F_{k}^{(i)}}{p_{n, n+1}}.
\de
Set
$$
\sigma_1=\sup_{n\geq 0}\sum_{k=0}^{n}
\frac{1}{p_{k, k+1}F_{k}^{(0)}}\sum_{j=n}^{\ift}F_{j}^{(0)}.$$

\begin{thm}\lb{miw}
If $\sigma_1<\ift$, then the skip-free Markov chain is geometrically transient.
\end{thm}
\bg{proof}
By Theorem \ref{zzz}, we need to construct a
solution to $(\ref{34e})$ for some $\lmd$, $b\in(0, 1)$. Let
\be\lb{lop1}
f_{i}=\l[p_{0 1}^{-1}\sum_{j=i}^{\ift}F_{j}^{(0)}\r]^{1/2}\q\mbox{and}\q
g_{i}=\sum_{j=i}^{\ift}F_{j}^{(0)}\sum_{k=0}^{j}
\frac{f_{k}}{p_{k, k+1}F_{k}^{(0)}},\q i\geq0.
\de
It is obvious that both $f$ and $g$ are decreasing. For $i\geq0$, define two operators
$$I_{i}(f)=\frac{F_{i}^{(0)}}{f_{i}-f_{i+1}}\sum_{k=0}^{i}\frac{f_{k}}
{p_{k, k+1}F_{k}^{(0)}}\q\mbox{and}\q
I\!I_{i}(f)=\frac{1}{f_{i}}\sum_{j=i}^{\ift}F_{j}^{(0)}
\sum_{k=0}^{j}\frac{f_{k}}{p_{k, k+1}F_{k}^{(0)}}.$$
Then by using the proportional property and \ct[Theorem 3.1]{cmf3}, we get
$$\sup_{i\geq0}I\!I_{i}(f)\leq \sup_{i\geq0}I_{i}(f)\leq 4\sigma_1.$$
Thus, combining \rf{lop1} with this inequality, we have $$\sup_{i\geq0}\frac{g_i}{f_i}=\sup_{i\geq0}I\!I_{i}(f)
\leq 4\sigma_1,$$ and
$$g_{0}=f_{0}I\!I_{0}(f)\leq f_{0}\sup_{i\geq0}I\!I_{i}(f)\leq 4\sigma_1 f_{0}=4\sigma_1\l[p_{0 1}^{-1}
\sum_{j=0}^{\ift}F_{j}^{(0)}\r]^{1/2}\leq 4\sgm_1^{3/2}.$$
We now determine $\lmd$, $b$ and a solution to inequality $(\ref{34e})$.
Set $\tilde{g}=g/g_{0}$. Then
\be\lb{plf1}\aligned
P \tld g(0)&=g_0^{-1}(p_{00}g_{0}+p_{01}g_{1})=1-p_{01}(g_0-g_{1})g_0^{-1}\\
&=1-f_0 g_0^{-1}\leq 1-\inf_{i\geq 0}\fr{f_i}{g_i}\leq 1-\fr{1}{4\sgm_1},\\
\endaligned
\de
and for $i\geq1$,
\be\lb{plf2}\aligned
P \tld g(i)&=g_0^{-1}\sum_{j=0}^{i+1}p_{i j}g_{j}
=g_0^{-1}\l[\sum_{k=0}^{i-1}\sum_{j=0}^{k}p_{i j}(g_{k}-g_{k+1})+p_{i, i+1}
g_{i+1}-p_{i, i+1}g_{i}+g_{i}\r]\\
&=g_0^{-1}\sum_{k=0}^{i-1}\sum_{j=0}^{k}p_{i j}F_{k}^{(0)}
\sum_{j=0}^{k}\frac{f_{j}}{p_{j, j+1}F_{j}^{(0)}}-g_0^{-1}p_{i, i+1}F_{i}^{(0)}
\sum_{j=0}^{i}\frac{f_{j}}{p_{j, j+1}F_{j}^{(0)}}+g_0^{-1}g_i\\
&\leq g_0^{-1}\sum_{k=0}^{i-1}p_{i}^{(k)}F_{k}^{(0)}\sum_{j=0}^{i-1}\frac{f_{j}}
{p_{j, j+1}F_{j}^{(0)}}-g_0^{-1}p_{i, i+1}F_{i}^{(0)}\sum_{j=0}^{i}\frac{f_{j}}
{p_{j, j+1}F_{j}^{(0)}}+g_0^{-1}g_i\\
&=g_0^{-1}p_{i, i+1}F_{i}^{(0)}\sum_{j=0}^{i-1}\frac{f_{j}}{p_{j, j+1}
F_{j}^{(0)}}-
g_0^{-1}p_{i, i+1}F_{i}^{(0)}\sum_{j=0}^{i}\frac{f_{j}}{p_{j, j+1}F_{j}^{(0)}}+g_0^{-1}g_i\\
&=\fr{g_i-f_i}{g_0}=\tld g_i-\fr{f_i}{g_i}\tld g_i\leq\tilde{g}_{i}-
\inf_{i\geq0}\frac{f_i}{g_i}~\tilde{g}_{i}\leq
\l(1-\frac{1}{4\sigma_1}\r)\tilde{g}_{i}.\\
\endaligned
\de
Therefore, combining \rf{plf1} with \rf{plf2},
$\tilde{g}$ is the nonnegative
solution of inequality $(\ref{34e})$ with $\lmd=b=1-\fr{1}{4\sgm_1}$. Hence the desired assertion follows.
\end{proof}

Next, we will study the algebraic transience by the first return time criteria \rf{cy4s}. For $i\geq0$ and $\ell\geq1$, let
$\tau_{0}=\inf\{n\geq1: \Phi_{n}=0\}$ and
$f_{i 0}^{(n)}=\P_i\{\tau_0=n\}$.
Set
$$m_{i 0}^{(0)}=\sum_{n=1}^{\ift}f_{i 0}^{(n)},\q
m_{i 0}^{(\ell)}=\sum_{n=1}^{\ift}n(n+1)\cdots(n+\ell-1)f_{i 0}^{(n)},$$
$$d_0^{(\ell)}=0,\q d_i^{(\ell)}=\sum_{k=1}^{i}\fr{F_i^{(k)}m_{k 0}^{(\ell-1)}}{p_{k, k+1}},\q
d^{(\ell)}=\sup_{i\geq1}\fr{\sum_{j=0}^{i-1}d_j^{(\ell)}}
{\sum_{j=0}^{i-1}F_j^{(0)}},$$
where $F_i^{(k)}$ is defined in \rf{000}. Set
$$\xi=\sup_{i\geq2}\fr{\sum_{j=1}^{i-1}F_j^{(0)}}{\sum_{j=0}^{i-1}F_j^{(0)}},\q
\sgm_2=\ell d^{(\ell)}.$$

\begin{thm}\lb{mzj1}
For the skip-free Markov chain, we have for $i\geq1$ and $\ell\geq1$,
$$m_{i 0}^{(0)}=\sum_{j=0}^{i-1}F_j^{(0)}\xi-\sum_{j=1}^{i-1}F_j^{(0)},$$
$$m_{i 0}^{(\ell)}=\ell\sum_{j=0}^{i-1}\l(F_j^{(0)}d^{(\ell)}-d_j^{(\ell)}\r),$$
and
$$m_{0 0}^{(0)}=p_{0 1}\xi+p_{0 0},\q m_{0 0}^{(\ell)}=p_{0 1}\sgm_2+\ell m_{0 0}^{(\ell-1)}.$$
Moreover, the chain is transient if and only if $\xi<1$; the chain is $\ell$-transient
if and only if $\xi<1$ and $\sgm_2<\ift$.
\end{thm}
\bg{rem}
By the Stolz theorem, it is obvious that $\xi<1$ if and only if $\sum_{j=0}^{\ift}F_{j}^{(0)}<\ift$, which is equivalent to
the transience by \ct[Theorem 4.52]{c1}.
\end{rem}
\bg{proof}
(1) Consider the following equations:
\be\lb{c1}
x_0=0,\q \sum_{k\not=0}p_{j k}x_k=x_j-p_{j 0},\q j\geq1.
\de
By the second successive
approximation for the minimal nonnegative solution,
\bg{equation}\lb{mta}
x_0=0\q\mbox{and}\q x_j=m_{j 0}^{(0)},\q j\geq1,
\end{equation}
is the minimal nonnegative solution of \rf{c1}.
Since $\sum_{k=0}^{j+1}p_{j k}=1$, and by the induction argument
$F_{j}^{(0)}=\sum_{k=1}^{j}p_k^{(0)}F_{j}^{(k)}/p_{k, k+1}$,
\rf{c1} can be rewritten as
\be\lb{pl5}\aligned
x_{j+1}-x_{j}&=\fr{1}{p_{j, j+1}}\l(\sum_{m=0}^{j-1}p_j^{(m)}(x_{m+1}-x_{m})-p_{j 0}\r)\\
&=\sum_{m=0}^{j-1}\fr{F_j^{(j)}p_j^{(m)}}{p_{j, j+1}}(x_{m+1}-x_{m})
-\fr{F_j^{(j)}p_{j 0}}{p_{j, j+1}}\\
&=\sum_{m=0}^{j-2}\fr{F_j^{(j)}p_j^{(m)}}{p_{j, j+1}}(x_{m+1}-x_{m})
+F_{j}^{(j-1)}(x_j-x_{j-1})
-\fr{F_j^{(j)}p_{j 0}}{p_{j, j+1}}\\
&=\sum_{m=0}^{j-2}\sum_{k=j-1}^{j}\fr{F_j^{(k)}p_k^{(m)}}{p_{k, k+1}}(x_{m+1}-x_{m})
-\sum_{k=j-1}^{j}\fr{F_j^{(k)}p_{k 0}}{p_{k, k+1}}=\cdots\\
&=\sum_{k=1}^{j}\fr{F_j^{(k)}p_k^{(0)}}{p_{k, k+1}}x_{1}
-\sum_{k=1}^{j}\fr{F_j^{(k)}p_{k 0}}{p_{k, k+1}}\\
&=F_j^{(0)}x_1-F_j^{(0)},\q j\geq 1.\nnb\\
\endaligned
\de
For $i\geq2$,
summing $j$ from 1 to $i-1$ gives
\be\lb{zst6}\aligned
x_i&=\sum_{j=0}^{i-1}F_j^{(0)}x_1-
\sum_{j=1}^{i-1}F_j^{(0)},\q i\geq2.
\endaligned
\de
Since \rf{mta} is the nonnegative solution of \rf{c1}, according to \rf{zst6}, we have
$m_{1 0}^{(0)}\geq\xi$.

On the other hand, set
$$u_0=0, \q u_1=\xi\q\mbox{and}\q u_i=\sum_{j=0}^{i-1}F_j^{(0)}\xi-\sum_{j=1}^{i-1}F_j^{(0)},\q i\geq 2.$$
Then $\{u_i\}$ is the nonnegative solution of \rf{c1}.
By the minimality of $m_{i 0}^{(0)}$, we get $m_{1 0}^{(0)}\leq\xi$.
Therefore, $m_{1 0}^{(0)}=\xi$, $m_{i 0}^{(0)}=u_i$ for $i\geq 2$, and
\bg{equation}\lb{b6w}
m_{0 0}^{(0)}=\sum_{i\not=0}p_{0 i}m_{i 0}^{(0)}+p_{0 0}=p_{0 1}\xi+p_{0 0}.\nnb
\end{equation}
Thus,
$\xi<1$ if and only if $m_{0 0}^{(0)}<1$,
which is equivalent to transience by Proposition \ref{ti8}.

(2) Consider the following equations:
\be\lb{c11}
x_0^{(\ell)}=0,\q \sum_{k\not=0}p_{j k}x_k^{(\ell)}
=x_j^{(\ell)}-\ell x_j^{(\ell-1)},\q \ell,\; j\geq1,\nnb
\de
where $x_j^{(0)}=m_{j 0}^{(0)}$. Similarly, $m_{1 0}^{(\ell)}=\sgm_2$,
$m_{i 0}^{(\ell)}=\ell\sum_{j=0}^{i-1}\l(F_j^{(0)}d^{(\ell)}-d_j^{(\ell)}\r)$
for $i\geq 2$, and
\bg{equation}\lb{b6w2}
m_{0 0}^{(\ell)}=\sum_{i\not=0}p_{0 i}m_{i 0}^{(\ell)}+\ell m_{0 0}^{(\ell-1)}=p_{0 1}\sgm_2+\ell m_{0 0}^{(\ell-1)}.\nnb
\end{equation}
Moreover,
by Corollary \ref{rmll}(2),
$$\E_0\l[\tau_0^{\ell}1_{\{\tau_0<\ift\}}\r]=\sum_{i\not=0}p_{0 i}\E_i\l[(\tau_0+1)^{\ell}1_{\{\tau_0<\ift\}}\r]+p_{00}=p_{0 1}\E_1\l[(\tau_0+1)^{\ell}1_{\{\tau_0<\ift\}}\r]+p_{00}.$$
From this, and
noting that there exist constants $c_1$ and $c_2$ such that
\bg{equation}\lb{mb9}
c_1 m_{1 0}^{(\ell)}\leq
\E_1\l[(\tau_0+1)^{\ell}1_{\{\tau_0<\ift\}}\r]\leq c_2 m_{1 0}^{(\ell)},\nnb
\end{equation}
we have $\E_0\l[\tau_0^{\ell}1_{\{\tau_0<\ift\}}\r]<\ift$ if and only if $m_{1 0}^{(\ell)}=\sgm_2<\ift$, which yields
the desired assertion by Theorem \ref{aaq}(2).
\end{proof}

Finally, we will study the criteria for strongly geometric transience and uniformly geometric transience of the sub-stochastic skip-free chain.
\begin{exm}\lb{zfu}
(The sub-stochastic skip-free chain).
Let $P=(p_{i j})$
be an irreducible sub-stochastic transition matrix on
$X=\N$ with
$p_{i j}=0$ for $j-i\geq2$, and
$\sum_{j\geq1}p_{i j}<1$ for some $i\geq1$.
\end{exm}
Let $\hat P=(\hat p_{i j})$ be a stochastic transition matrix on the state space $\hat X=X\cup\{0\}$,
where
\be\lb{s35}\aligned\widehat{p}_{i j}=
    \left\{\begin{array}{ll}
     p_{i j},    & i, j\in X;\\
     1-\sum_{j\geq1}p_{i j},  & i\in X,\, j=0;\\
     1,          & i, j=0;\\
     0,          & i=0,\, j\in X.\nnb\\
   \end{array}    \right.\\
   \endaligned
   \de
Then $0$ is an absorbing state for $\hat P$. Thus, $\tau$ defined in \rf{876}
is just the first hitting time of state $0$ for $\hat P$, say $\hat\tau_0$.
For $i\geq1$, define
$$d_0=0,\q d_i=\sum_{k=1}^{i}\fr{F_i^{(k)}}{p_{k, k+1}}\q\mbox{and}\q
d=\sup_{i\geq1}\fr{\sum_{j=0}^{i-1}d_j}{\sum_{j=0}^{i-1}F_j^{(0)}}.$$
Let
$$\sgm_3=\sup_{n\geq1}\sum_{k=0}^{n-1}F_{k}^{(0)}\sum_{j=n}^{\ift}\fr{1}{p_{j, j+1}F_{j}^{(0)}}
\q\mbox{and}\q\sgm_4=\sup_{n\geq0}\sum_{k=0}^{n}\l(F_k^{(0)}d-d_k\r).$$
Following the argument of single-birth processes in \ct{zyh1} or \ct[Theorem 4.52]{c1}, we can derive that there exists a constant $\kp>1$ such that for all $i\geq1$,
$$\E_i\kp^{\tau}=\E_i\kp^{\hat\tau_0}<\ift$$
provided $\sgm_3<\ift$,
which also implies that $$\P_i\l\{\tau<\ift\r\}=\P_i\l\{\hat\tau_0<\ift\r\}=1,\q
i\geq1.$$
Meanwhile, $$\sup_{i\geq1}\E_i\tau=\sup_{i\geq1}\E_i\hat\tau_0<\ift$$
if and only if $\sgm_4<\ift$.
Therefore, we have the following criteria according to Theorems \ref{28}(2) and \ref{39}(5).

\bg{thm}\lb{bwo}
$(1)$ The chain $\Phi$ is strongly geometric transient if $\sgm_3<\ift$.

$(2)$ The chain $\Phi$ is uniformly geometric transient if and only if $\sgm_4<\ift$.
\end{thm}

\bigskip
\bigskip
\noindent\textbf{Acknowledgements}
The authors would thank Professor Mu-Fa Chen for introducing us the topic on
transience for Markov chains, and thank Professor Yu-Hui Zhang for his help on the skip-free chains.
This work is supported in part
by 985 Project, 973 Project (No 2011CB808000), NSFC (No 11131003),
SRFDP (No 20100003110005) and the Fundamental Research Funds for the Central Universities.

\end{document}